\documentclass[11pt,oneside,british]{article}
\newcommand{\comment}[1]		{}
\usepackage{amsmath,amssymb,amsthm}
\usepackage{rotating,pdflscape,afterpage}
\usepackage[sc,osf]{mathpazo}
	\AtBeginDocument{
		\DeclareSymbolFont{AMSb}{U}{msb}{m}{n}
		\DeclareSymbolFontAlphabet{\mathbb}{AMSb}
	}
\usepackage{graphicx}
\usepackage{extpfeil}
\usepackage{mathabx} %UNDERBRACE BREAKS IN THE OPPOSITE ORDER
\usepackage{fancyvrb}
\usepackage{float}
\usepackage{caption}
\newcommand{\mockalph}[1]{\!}
\usepackage{enumitem}
\usepackage[nouppercase]{scrlayer-scrpage}
\usepackage{avant} 
\usepackage{mathrsfs}
\usepackage{mathtools}
\usepackage{faktor,xfrac}
\usepackage[all]{xypic}\xyoption{rotate}
	\newdir{ >}{{}*!/-7pt/\dir{>}}
	\newdir{i}{{}*!/-7pt/\dir{(}	}
\usepackage{thmtools}

\usepackage[utf8]{inputenc}
\usepackage[T1]{fontenc}
\usepackage{tocloft}
\makeatletter
\renewcommand{\l@figure}{\@dottedtocline{1}{1em}{3.5em}}
\renewcommand{\l@table}{\@dottedtocline{2}{1em}{3.5em}}
\makeatother
\newcommand*{\noaddvspace}{\renewcommand*{\addvspace}[1]{}}
\addtocontents{lof}{\protect\noaddvspace}

\usepackage[hyphens]{url}		%
\usepackage{hyperref}	%EVERYTHING BREAKS IF THESE THREE ARE >
\usepackage{cleveref} 	%PUT IN ANOTHER ORDER.

\makeatletter
\let\c@figure\c@table
\let\c@equation\c@table
\makeatother

\numberwithin{table}{section}
\numberwithin{figure}{section}
\newtheorem{abcthm}[table]{Theorem}

\newtheorem{theorem}[table]{Theorem}
\newtheorem{thmx}{Theorem}

\newtheorem{proposition}[table]{Proposition}

\newtheorem{corollary}[table]{Corollary}

\newtheorem{lemma}[table]{Lemma}

\newtheorem{claim}[table]{Claim}

\theoremstyle{definition}

\newtheorem{definition}[table]{Definition}

\newtheorem{notation}[table]{Notation}

\newtheorem{observation}[table]{Observation}

\newtheorem{conjecture}[table]{Conjecture}

\theoremstyle{remark}
\newtheorem{fact}[table]{Fact}

\newtheorem{example}[table]{Example}
\newtheorem{exercise}[table]{Exercise}
\newtheorem{EG}[table]{Example}

\newtheorem{problem}[table]{Problem}
\newtheorem{histrmks}[table]{Historical remarks}
\newtheorem{remark}[table]{Remark}
\newtheorem{remarks}[table]{Remarks}

\theoremstyle{plain}
\newtheorem*{thm*}{Theorem}
\newtheorem*{theorem*}{Theorem}
\newtheorem*{prop*}{Proposition}
\newtheorem*{proposition*}{Proposition}
\newtheorem*{lemma*}{Lemma}
\newtheorem*{corollary*}{Corollary}
\newtheorem*{cor*}{Corollary}

\theoremstyle{definition}
\newtheorem*{definition*}{Definition}
\newtheorem*{defn*}{Definition}
\newtheorem*{QQ*}{Question}
\newtheorem*{obs*}{Observation}
\newtheorem*{notation*}{Notation}

\theoremstyle{remark}
\newtheorem*{rmk*}{Remark}
\newtheorem*{remark*}{Remark}

\newtheorem*{examples*}{Examples}
\newtheorem*{example*}{Example}
\newtheorem*{EG*}{Example}
\newtheorem*{EGs*}{Examples}
\newtheorem*{fact*}{Fact}
\newtheorem*{prob*}{Problem}

\usepackage{etoolbox}
\usepackage[usenames,dvipsnames,svgnames,table]{xcolor}
\newcommand		{\defd}[1]	{\textcolor{RoyalBlue}{\textbf{\textit{#1}}}}
\newcommand		{\defm}[1]	{\textcolor{RoyalBlue}{#1}}

\tracingpatches
\makeatletter
\patchcmd{\@setref}{\bfseries ??}{\bfseries\color{red} FIX ME!}{}{}
\patchcmd{\@setcite}{\bfseries ?}{\bfseries\color{red} FIX ME!}{}{}
\patchcmd{\@setcref}         {??}{\color{red} FIX ME!}{}{}
\patchcmd{\@setcref}         {??}{\color{red} FIX ME!}{}{}
\patchcmd{\@setcrefrange}    {??}{\color{red} FIX ME!}{}{}
\patchcmd{\@setcrefrange}    {??}{\color{red} FIX ME!}{}{}
\patchcmd{\@setcrefrange}    {??}{\color{red} FIX ME!}{}{}
\patchcmd{\@setcrefrange}    {??}{\color{red} FIX ME!}{}{}
\patchcmd{\@setcrefrange}    {??}{\color{red} FIX ME!}{}{}
\patchcmd{\@setcrefrange}    {??}{\color{red} FIX ME!}{}{}
\patchcmd{\@setnamecref}     {??}{\color{red} FIX ME!}{}{}
\patchcmd{\@setnamecref}     {??}{\color{red} FIX ME!}{}{}
\patchcmd{\@setcpageref}     {??}{\color{red} FIX ME!}{}{}
\patchcmd{\@setcpageref}     {??}{\color{red} FIX ME!}{}{}
\patchcmd{\@setcpagerefrange}{??}{\color{red} FIX ME!}{}{}
\patchcmd{\@setcpagerefrange}{??}{\color{red} FIX ME!}{}{}
\patchcmd{\@setcpagerefrange}{??}{\color{red} FIX ME!}{}{}
\patchcmd{\@setcpagerefrange}{??}{\color{red} FIX ME!}{}{}
\patchcmd{\@setcpagerefrange}{??}{\color{red} FIX ME!}{}{}
\patchcmd{\@cref}            {??}{\color{red} FIX ME!}{}{}
\def\blx@citation@entry#1#2{%
  \blx@bibreq{#1}%
  \ifinlist{#1}{\blx@cites}
    {}
    {\listgadd{\blx@cites}{#1}%
     \blx@auxwrite\@mainaux{}{\string\abx@aux@cite{#1}}}%
  \ifinlistcs{#1}{blx@segm@\the\c@refsection @\the\c@refsegment}
    {}
    {\listcsgadd{blx@segm@\the\c@refsection @\the\c@refsegment}{#1}}%
  \blx@ifdata{#1}%
    {}%
    {\ifcsdef{blx@miss@\the\c@refsection}%
       {\ifinlistcs{#1}{blx@miss@\the\c@refsection}%
          {{\bfseries\color{red} cite:} }%
          {\blx@logreq@active{#2{#1}}}}%
       {\blx@logreq@active{#2{#1}}}}}
\def\blx@citeadd#1{%
  \ifcsdef{blx@keyalias@\the\c@refsection @#1}
    {\edef\blx@realkey{\csuse{blx@keyalias@\the\c@refsection @#1}}}
    {\def\blx@realkey{#1}}%
  \expandafter\blx@citation\expandafter{\blx@realkey}\blx@msg@cundefon
  \expandafter\blx@ifdata\expandafter{\blx@realkey}
    {\advance\blx@tempcnta\@ne
     \listeadd\blx@tempa{\blx@realkey}}
    {\ifnum\blx@tempcntb>\z@\multicitedelim\fi
     \expandafter\abx@missing\expandafter{\blx@realkey}%
     \advance\blx@tempcntb\@ne}}

\DeclarePairedDelimiterX{\pmodx}[1]{(}{)}{{\operator@font mod}\mkern6mu#1}
\renewcommand{\pmod}{%
  \allowbreak
  \if@display\mkern18mu\else\mkern8mu\fi
  \pmodx
}
\makeatother

\makeatletter
\newcommand{\oset}[3][0ex]{%
\raisebox{.175ex}{$%
  \mathrel{\mathop{#3}\limits^{
    \vbox to#1{\kern-2\ex@
    \hbox{$\scriptstyle#2$}\vss}}}
    $}%
    }
\makeatother

\newcommand{\myred}{BrickRed}
\hypersetup{
		colorlinks			= true, 	
		urlcolor			= \myred, 
		linkcolor			= \myred, 	
		citecolor			= PineGreen 	
}

\usepackage{xspace}

\usepackage[left=2.54cm,right=2.54cm,top=2.54cm,bottom=2.54cm]{geometry}
\usepackage{tikz}
\usetikzlibrary{matrix,fit,backgrounds,calc,arrows} %ARXIV HATES arrows.meta
\tikzstyle{image}=[rectangle,fill=Red!20,inner sep=-2pt]
\tikzstyle{nonzero}=[rectangle,fill=Navy!20,inner sep=0pt]
\tikzstyle{nonzerosm}=[rectangle,fill=Navy!20,inner sep=-2pt]

\makeatletter
\newbox\xrat@below
\newbox\xrat@above
\newcommand{\xrightarrowtail}[2][]{%
  \setbox\xrat@below=\hbox{\ensuremath{\scriptstyle #1}}%
  \setbox\xrat@above=\hbox{\ensuremath{\scriptstyle #2}}%
  \pgfmathsetlengthmacro{\xrat@len}{max(\wd\xrat@below,\wd\xrat@above)+.6em}%
  \mathrel{\tikz [>->,baseline=-.55ex]
                 \draw (0,0) -- node[below=-2pt] {\box\xrat@below}
                                node[above=-2pt] {\box\xrat@above}
                       (\xrat@len,0) ;}}
\newbox\xrat@below
\newbox\xrat@above
\renewcommand{\xtwoheadrightarrow}[2][]{%
  \setbox\xrat@below=\hbox{\ensuremath{\scriptstyle #1}}%
  \setbox\xrat@above=\hbox{\ensuremath{\scriptstyle #2}}%
  \pgfmathsetlengthmacro{\xrat@len}{max(\wd\xrat@below,\wd\xrat@above)+.6em}%
  \mathrel{\tikz [->>,baseline=-.55ex]
                 \draw (0,0) -- node[below=-2pt] {\box\xrat@below}
                                node[above=-2pt] {\box\xrat@above}
                       (\xrat@len,0) ;}}
\makeatother
\newcommand{\xmono}{\xrightarrowtail}
\newcommand{\mono}{\xmono{\phantom{\ \, }}}
\newcommand{\xepi}{\xtwoheadrightarrow}

\makeatletter
\newcommand{\presectionskip}{-1.5\baselineskip}
\newcommand{\postsectionskip}{0.3\baselineskip}
\usepackage{titlesec}
\renewcommand{\section}{\@startsection
  {chapter}{0}{0mm}%                  % name, level, indent
  {\presectionskip}%                   % beforeskip
  {\postsectionskip}%                 % afterskip
  {\sffamily\huge}}% % style		\HUGE ?%\bfseries
\renewcommand{\section}{\@startsection
  {section}{1}{0mm}%                  % name, level, indent
  {\presectionskip}%                   % beforeskip
  {\postsectionskip}%                 % afterskip
  {\sffamily\LARGE}}% % style%\bfseries
\renewcommand{\subsection}{\@startsection
  {subsection}{2}{0mm}%                  % name, level, indent
  {\presectionskip}%                   % beforeskip
  {\postsectionskip}%                 % afterskip
  {\sffamily\Large}}% % style%\bfseries
\renewcommand{\subsubsection}{\@startsection
  {subsubsection}{3}{0mm}%                  % name, level, indent
  {\presectionskip}%                   % beforeskip
  {\postsectionskip}%                 % afterskip
  {\sffamily\normalsize}}% % style%\bfseries
\renewcommand{\@seccntformat}[1]{\csname the#1\endcsname.\quad}
\newcommand\HUGE{\@setfontsize\Huge{30}{47}} 
\titleformat{\chapter}[display]
{\sffamily\Large}
{Chapter {\HUGE\normalfont\thechapter}}    %% change here
{1em}
{\huge}
\titleformat{\section}
{\sffamily\Large}
{\normalfont{\@setfontsize\Large{18}{47}\thesection.}}    %% change here
{1em}
{}
\titleformat{\subsection}
{\sffamily\large}
{\normalfont{\Large\thesubsection}.}    %% change here
{1em}
{}
\titleformat{\subsubsection}
{\sffamily\normalsize}
{\normalfont{\large\thesubsubsection}.}    %% change here
{1em}
{}
\makeatother

\makeatletter
\def\smallunderbrace#1{\mathop{\vtop{\m@th\ialign{##\crcr
   $\hfil\displaystyle{#1}\hfil$\crcr
   \noalign{\kern3\p@\nointerlineskip}%
   \tiny\upbracefill\crcr\noalign{\kern3\p@}}}}\limits}
\makeatother

\renewcommand{\SS}{\textsection}
\newcommand{\bthm}{\begin{theorem}}
\newcommand{\ethm}{\end{theorem}}
\newcommand{\bprop}{\begin{proposition}}
\newcommand{\eprop}{\end{proposition}}
\newcommand{\bcor}{\begin{corollary}}
\newcommand{\ecor}{\end{corollary}}
\newcommand{\bconj}{\begin{conjecture}}
\newcommand{\econj}{\end{conjecture}}
\newcommand{\blem}{\begin{lemma}}
\newcommand{\elem}{\end{lemma}}
\newcommand{\bclm}{\begin{claim}}
\newcommand{\eclm}{\end{claim}}
\newcommand{\bpf}{\begin{proof}}
\newcommand{\epf}{\end{proof}}
\newcommand{\bdetails}{\begin{details}}
\newcommand{\edetails}{\end{details}}
\newcommand{\bdefi}{\begin{definition}}
\newcommand{\edefi}{\end{definition}}
\newcommand{\bdefn}{\begin{definition}}
\newcommand{\edefn}{\end{definition}}
\newcommand{\bex}{\begin{example}}
\newcommand{\eex}{\end{example}}
\newcommand{\bprob}{\begin{problem}}
\newcommand{\eprob}{\end{problem}}
\newcommand{\bob}{\begin{observation}}
\newcommand{\eob}{\end{observation}}
\newcommand{\bexer}{\begin{exercise}}
\newcommand{\eexer}{\end{exercise}}
\newcommand{\bexers}{\begin{exercises}}
\newcommand{\eexers}{\end{exercises}}
\newcommand{\brmk}{\begin{remark}}
\newcommand{\ermk}{\end{remark}}
\newcommand{\bhist}{\begin{histrmks}}
\newcommand{\ehist}{\end{histrmks}}
\newcommand{\brmks}{\begin{remarks}}
\newcommand{\ermks}{\end{remarks}}
\newcommand{\bverb}{\begin{verbatim}}
\newcommand{\everb}{\end{verbatim}}
\newcommand{\bntn}{\begin{notation}}
\newcommand{\entn}{\end{notation}}
\newcommand{\bfct}{\begin{fact}}
\newcommand{\efct}{\end{fact}}
\newcommand{\bfcts}{\begin{facts}}
\newcommand{\befcts}{\end{facts}}
\newcommand{\benum}{\begin{enumerate}}
\newcommand{\eenum}{\end{enumerate}}
\newcommand{\bitem}{\begin{itemize}}
\newcommand{\eitem}{\end{itemize}}

\renewcommand	{\epsilon}	{\varepsilon}
\renewcommand	{\a}		{\alpha}
\renewcommand	{\b}		{\beta}
\renewcommand	{\d}		{\delta}
\newcommand		{\h}		{\eta}

\newcommand		{\vk}		{\varkappa}

\newcommand		{\gpp}[1]	{\mn\big( #1 \big)}

\renewcommand	{\:}		{\colon}

\newcommand		{\quotientmed}[2]	{{\raisebox{.2em}{$#1$}}\  \!\!\big/\!\!\ 
										{\raisebox{-.2em}{$#2$}}}

\newcommand		{\fk}		{{\mathfrak k}}

\newcommand		{\fg}		{{\mathfrak g}}

\newcommand		{\td}		{\tilde{d}}

\newcommand		{\wP}		{\widehat{P}}
\newcommand		{\cP}		{\widecheck{P}}
\newcommand		{\ewP}	{\ext\wP}

\newcommand		{\eP}		{\ext P}

\newcommand		{\ecP}	{\ext\cP}

\newcommand		{\LSS}	{Leray spectral sequence\xspace}
\newcommand		{\SSS}	{Serre spectral sequence\xspace}
\newcommand		{\EMSS}	{Eilenberg--Moore spectral sequence\xspace}

\newcommand		{\exterior}	{\Lambda}
\newcommand		{\ext}		{\exterior}

\newcommand		{\ii}		{\iota}
\newcommand		{\susp}		{\Sigma}

\renewcommand	{\th}			{^{\mathrm{th}}}

\newcommand		{\DGA}		{\textsc{dga}\xspace}

\newcommand		{\CDGA}		{\textsc{cdga}\xspace}

\newcommand{\eQ}{\ext Q}

\newcommand		{\dual}	{\mn^\vee}

\newcommand		{\APL}	{A_{\mathrm{PL}}}

\makeatletter
\newcommand{\subalign}[1]{%
  \vcenter{%
    \Let@ \restore@math@cr \default@tag
    \baselineskip\fontdimen10 \scriptfont\tw@
    \advance\baselineskip\fontdimen12 \scriptfont\tw@
    \lineskip\thr@@\fontdimen8 \scriptfont\thr@@
    \lineskiplimit\lineskip
    \ialign{\hfil$\m@th\scriptstyle##$&$\m@th\scriptstyle{}##$\crcr
      #1\crcr
    }%
  }
}
\makeatother

\newcommand		{\bl}		{\bullet}

\newcommand		{\floor}[1]			{\left\lfloor #1 \right\rfloor}

\newcommand		{\eqn}[1]			{\begin{align*} #1 \end{align*}}
\newcommand		{\quation}[1]		{\begin{equation} #1 \end{equation}}

\newcommand		{\case}[1]			{\begin{cases} #1 \end{cases}}

\newcommand		{\hyref}[1]			{\hyperref[#1]{\ref*{#1}}}

\newcommand		{\mn}				{\mspace{-2mu}}
\newcommand		{\mnn}				{\mspace{-1mu}}
\newcommand		{\dsp}			{\displaystyle}

\newcommand		{\nd}				{\noindent}

\newcommand		{\os}			{\overset}
\newcommand		{\us}			{\underset}

\newcommand		{\wt}			{\widetilde}

\newcommand		{\mr}			{\mathrm}
\newcommand		{\bb}			{\mathbb}

\newcommand		{\f}			{\mathfrak}
\newcommand		{\ms}			{\mathscr}

\newcommand		{\g}			{\gamma}
\newcommand		{\e}			{\epsilon}
\newcommand		{\z}			{\zeta}

\renewcommand		{\L}			{\ms{L}}
\newcommand		{\s}			{\sigma}

\newcommand		{\G}			{\Gamma}
\newcommand		{\D}			{\Delta}

\DeclareSymbolFont{cmletters}{OT1}{cmr}{m}{n}
\DeclareMathSymbol{\Ups}{\mathalpha}{cmletters}{"7}
\renewcommand	{\Upsilon}{\Ups}

\newcommand		{\Z}		{\bb Z}
\newcommand		{\Q}		{\bb Q}
\newcommand		{\R}		{\bb R}

\DeclareMathOperator{\id}			{id}
\DeclareMathOperator{\pr}		{pr}

\DeclareMathOperator{\Sq}{Sq}

\renewcommand 		{\H}		{H^*}
\newcommand 		{\HG}	{{\H_G}}

\newcommand 		{\HT}	{{\H_T}}

\newcommand 		{\HK}	{{\H_K}}

\newcommand		{\KGK}	{{}_K G_K}

\DeclareRobustCommand{\lq}		{\text{\reflectbox{$/$}}}		%LEFT QUOTIENT

\renewcommand	{\-}		{^{-1}}

\renewcommand	{\.}		{\cdot}
\newcommand		{\x}		{\times}

\newcommand		{\semidirect}	{\rtimes}
\newcommand{\oplushigher}{\mathbin{\raisebox{.85pt}{$\displaystyle\oplus$}}}
\newcommand{\otimeshigher}{\mathbin{\raisebox{.85pt}{$\displaystyle\otimes$}}}
\DeclareMathOperator*{\otimesvariable}{%
			\mathchoice {\raisebox{.85pt}{$\displaystyle\otimes$}}
						{\raisebox{.85pt}{$\otimes$}}
						{\raisebox{0.7pt}{$\scriptstyle\otimes$}}
						{\raisebox{0.2pt}{$\scriptscriptstyle\otimes$}}
						}
\newcommand		{\tensor}		{\otimesvariable}
\newcommand		{\xt}[3]		{{#2}\us{#1}\otimeshigher{#3}}
\newcommand		{\direct}		{\oplushigher}
\newcommand		{\ox}			{\tensor}

\newcommand		{\+}			{\direct}

\DeclareMathOperator{\diag}		{diag}

\DeclareMathOperator{\rk}			{rk }
\DeclareMathOperator{\im}		{im }

\DeclareMathOperator{\Tor}		{Tor}

\DeclareMathOperator{\Ad}		{Ad }

\DeclareMathOperator{\Hom}		{Hom}

\renewcommand	{\O}			{\mr{O}}
\newcommand		{\Orth}		{\O}
\newcommand		{\SO}		{\mr{SO}}

\newcommand		{\U}			{\mr{U}}
\newcommand		{\SU}			{\mr{SU}}

\newcommand		{\Sp}			{\mr{Sp}}
\newcommand		{\Spin}		{\mr{Spin}}

\newcommand		{\longto} 		{\longrightarrow}
\newcommand		{\lt}			{\longto}

\newcommand		{\lmt}			{\longmapsto}

\newcommand		{\inc}		{\hookrightarrow}
\newcommand		{\xinc}		{\xhookrightarrow}
\newcommand		{\longinc}		{\xinc[]{\ \ \ \ }}

\newcommand		{\longmono}	{\xmono[]{\ \ \ \ }}

\newcommand		{\longepi}	{\xepi[]{\ \ \ \ }}

\newcommand		{\simto}		{\xrightarrow{\sim}}

\newcommand		{\longsimto}	{\os\sim\longto}
\newcommand		{\isoto}		{\longsimto}

\newcommand		{\ceq}			{\coloneqq}

\newcommand		{\iso}				{\cong}
\newcommand		{\homeo}			{\approx}

\numberwithin{equation}{section}
\theoremstyle{definition}
\newtheorem{discussion}[table]{Discussion}
\newtheorem{history}[table]{History}
\newtheorem{construction}[table]{Construction}
\newcommand{\pont}{A}
\newcommand{\redpont}{B}
\newcommand{\Zt}{\Z[\sfrac 1{\,2}]}

\newcommand{\HKz}{\H_{K_0}}
\newcommand{\HKzt}{H^{\ox 2}_{K_0}}
\newcommand{\KzGK}{{}_{K_0}\mspace{-2mu} G_{K}}
\newcommand{\KzGKz}{{}_{K_0}\mspace{-2mu} G_{K_0}}
\newcommand{\forms}{A}
\newcommand{\thesubalg}{A'}
\newcommand{\bas}{_{\mathrm{bas}}}

\usepackage{tgpagella}

\title{Grassmannians and the equivariant cohomology of 
	isotropy actions
}
\author{Jeffrey D. Carlson}
\begin{document}
\maketitle

\begin{abstract}
We prove a general structure theorem
for the rational Borel equivariant cohomology ring
of an equivariantly formal isotropy action,
rederiving He's recent computation
of the equivariant cohomology of real Grassmannians 
as an illustration.
%the recent Recent work of He has determined the rational
%Borel s of the isotropy actions
%on real Grassmannians through GKM methods.
%n this note,
%we propound two less involved approaches,
%one obtaining the result as a corollary of other known theorems
%on the equivariant cohomology of isotropy actions on homogeneous spaces
%%due to Kumar--Vergne and the author,
%and the other reducing it to a calculation 
%in a more generally applicable 
%Particularly, this gives effortless calculations 
%of the equivariant cohomology of the isotropy actions on symmetric spaces.
\end{abstract}

\renewcommand{\thefootnote}{\fnsymbol{footnote}} 
\renewcommand{\thefootnote}{\arabic{footnote}}

%\section*{Introduction}
The rational singular cohomology rings of
real Grassmannians have long been
known,
%~\cite{lerayCR1949e,takeuchi1962pontrjagin}\cite[p.~71]{cartan1950transgression},
but surprisingly
%These can be viewed as homogeneous quotients of 
%the special orthogonal group
%\eqn{
%	\tG_\ell(\R^{m+\ell}) = G/K_0 &= \,
%		\quotientmed{\SO(m+\ell)\,}
%			{\,\big(\SO(m) \x \SO(\ell)\big)},\\
%	G_\ell(\R^{m+\ell})   = G/K   &= \,
%		\quotientmed{\SO(m+\ell)\,}
%			{\,\mathrm S\big(\O(m) \x \O(\ell)\big)},
%}
%which admit a left \emph{isotropy action}
%by the stabilizer of a point.
the rational Borel equivariant cohomology rings 
of the corresponding isotropy actions
%on oriented real Grassmannians 
were only recently published (by He,
then again by Sadykov)~\cite[5.2.1--2,\, 6.3.1]{he2016grassmannian}\cite[Thm.~6]{sadykov2017}.
This present note reproves this calculation
%as well as that for the unoriented counterparts,
as a consequence 
of a general structure theorem for such actions
obtainable through a Sullivan model due to V. Kapovitch~\cite{kapovitch2002biquotients}.

% example of a more general method, using 
%and which we take this opportunity to better propagandize.

\begin{thmx}
	\label{thm:cohomisotropyformal}
	Let $G$ be a connected Lie group and $H$ 
	a closed, connected subgroup 
	such that the left multiplication action of $H$ on $G/H$
	is equivariantly formal.
	Then there is a natural ring isomorphism 
	\[
	\H_H(G/H;\Q) \iso \xt{\H(BG;\Q)}{\H(BH;\Q)}{\H(BH;\Q)}\, \ox_\Q\, \ewP,
	\]
	where $\ewP$ is 
	%an exterior algebra
	isomorphic to the image of $\H(G/H;\Q) \lt \H (G;\Q)$.
	Particularly, $\dim_\Q \wP = \rk G - \rk H$. 
\end{thmx}

It is known that the left multiplication action of $H$ 
on the right quotient $G/H$,
the so-called \defd{isotropy action},
is equivariantly formal 
\benum
\item when $G$ and $H$ 
share a maximal torus~\cite[Thm.~14.1(1)]{GKM1998},
\item when $\H(G;\Q) \lt \H(H;\Q)$ 
is surjective~\cite[Thm.~A, Prop.~4.1]{shiga1996equivariant}, 
\item when $H \iso S^1$ and $N_G(H)$ has two components%
	~\cite[Prop.~1.4]{carlson2018eqftorus},%
\footnote{\ 
	Such an action is equivariantly formal if and only if
	this happens or $\H(G;\Q) \lt \H(H;\Q)$ is surjective,
	and these possibilities are exclusive.
	}
and
\item when $G/H$ is a 
\defd{generalized symmetric space}~\cite{goertschesnoshari2016}
in the sense that $H$ is the identity component of the fixed point set
of a smooth group automorphism of $G$.
\eenum
The first three of these cases lead to specialized structure theorems 
(\Cref{thm:equalrank,thm:Hsurj,thm:circle}),
one of which is standard by other methods~\cite[Prop.~68, p.~161]{kumarvergne1993equivariant}.
The He result follows from \Cref{thm:cohomisotropyformal}
as Grassmannians are symmetric spaces, and
the variants for unoriented Grassmannians follow by a covering space argument.
A certain amount of notation has to be established to
%in \Cref{rmk:bignotation} 
state these ring structures in closed
form in \Cref{thm:He}, and we delay this work 
and hence the statement of the result until \Cref{sec:Grass}.
	We also improve the coefficient ring to $\Zt$.
%
%But these structures in fact follow from a more general theorem,
%%This result appeared in an a superseded version of the present
%%author's dissertation--turned--book project, 
%which has
%and has been mentioned in published work in one other paper~\cite[]{CF...},
%but no proof has been published

As several of these results have been rediscovered
multiple times and could in principle have been known
since Borel's thesis,
we include an abundance of historical clarification
(\ref{rmk:Borel} and \ref{rmk:history}).

\section{Models}\label{sec:models}

\Cref{thm:cohomisotropyformal} 
is based on a rational model for
homotopy biquotients to be found in an unpublished manuscript of Kapovitch~\cite[Prop.~1]{kapovitch2002biquotients}
and later in a textbook of 
F\'{e}lix \emph{et al.}~\cite[Thm.~3.50]{FOT}\footnote{\ 
	and independently, much later, by the present author,
	lest it be thought that he is immune to 
	the reduplication of effort described in \Cref{rmk:history}
}
and modeled after results of 
Eschenburg~{\cite[Thm.~1, p.~159]{eschenburg1992biquotient}},
which we rederive briefly here.
All these models will be pure Sullivan algebras.

\bdefn
	In what follows,
\textbf{the coefficient ring is always $\Q$ unless otherwise specified}. 
The ring over which a tensor product is taken will always
be clear from context.
A \defd{pure Sullivan algebra} $(\eQ \ox \eP,d)$
is a finitely-generated commutative differential graded algebra
(\defd{\fontfamily{qpl}\selectfont \textsc{cdga}})
over the rationals,
with $\eP$ an exterior algebra on a 
graded rational vector space $P$ of odd-degree generators,
$\eQ$ a symmetric algebra on an evenly-graded space of generators,
and $d$ a derivation of degree $1$ such that 
\[
dP \leq \eQ \qquad \textrm{and}\qquad d\eQ = 0.
\]
A \defd{pure Sullivan model} of a space $X$ is a $\Q$-\CDGA map from
a pure Sullivan algebra %$(\eQ \ox \eP,d)$ 
to the algebra $\defm{\APL(X)}$
of polynomial differential forms on $X$.
All we need to know about the latter is that it is a $\Q$-\CDGA
functorial in spaces
which can be connected to the singular cochain algebra $C^*(X;\Q)$
by a zigzag of natural $\Q$-\DGA quasi-isomorphisms.
A space $X$ is \defd{formal} if there is a zigzag of $\Q$-\CDGA
quasi-isomorphisms between $\APL(X)$ 
and $(\H X,0)$.
In what follows we will abusively refer to models in terms of their source
algebras and take the maps to $\APL(X)$ for granted.
\edefn

Crucially, Sullivan models behave well with respect to fibrations and pullbacks.\footnote{\ 
	We only need the following result for the \emph{pure} Sullivan models
	we have defined, but it holds of all Sullivan models.					} %%%

\bthm[{\cite[Prop.~15.5,8]{FHT}}]\label{thm:FHTfibration}
Given a map of Serre fibrations
\[
\xymatrix{
	F\ar[r]\ar[d] 	& F'\ar[d]\\
	E\ar[r]\ar[d] 	& E'\ar[d]^q\\
	B\ar[r]_f		& B'
}
\]
and Sullivan models $(\ext V_{B'},d) \lt (\ext V_{B},d)$ for $f$
and $(\ext V_{B'},d) \lt (\ext V_{B'} \ox \ext V_{F'},d)$ for $q$,
if $\H F' \lt \H F$ is an isomorphism,
$\pi_1 B = \pi_1 B' = \pi_0 E = \pi_0 E' = 0$,
and either $\H F$
or both of $\H B$ and $\H B'$ are of finite type,
then $E$ admits a Sullivan model % $(\ext V_{E},d) \lt \APL(E)$ with
\[
(\ext V_{E},d)
=
\xt{(\ext V_{B'},d)}{(\ext V_{B},d)}{(\ext V_{B'} \ox \ext V_{F'},d)} 
\iso 
(\ext V_{B} \ox \ext V_{F'},d).
\] 
\ethm

We will apply this theorem to find a model for
the Borel equivariant cohomology of a translation action on a homogeneous space.

\begin{definition}\label{def:actions}
	If $G$ is a topological group,
	we will write $\defm{EG \to BG}$ for any universal principal $G$-bundle.
	When $G$ is a Lie group and $H$ a closed
	subgroup, then $EG \to EG/H$ is also a universal principal $H$-bundle.
	We allow the action of $G$ or a subgroup on $EG$ to be written
	on either side; 
	one can employ the standard trick $g\.e \ceq eg\-$
	or note that the standard join and simplicial space constructions of $EG$
	actually admit actions of $G$ on both sides.
	
	Given a continuous left action of $G$ on a space $X$,
	we write 
	\[
	\defm {X_G} \ceq EG \ox_G X = \frac{EG \x X}{(eg,x) \sim (e,gx)}.
	\]
	for the \defd{homotopy quotient} or \emph{homotopy orbit space}
	and $\defm{e \ox x}$ for its elements.
	Its singular cohomology $\defm{\HG X}$
	is the \defd{(Borel) equivariant cohomology} of the action.
	The Borel fibration $X \to X_G \to BG$,
	given by $x \lmt e_0 \ox x$ for any fixed $e_0 \in EG$
	followed by $e \ox x \lmt eG$,
	induces ring maps $\defm{\HG} \ceq \H(BG) \to \HG X \to \H X$
	in singular cohomology (implicitly 
	with respect to $\Q$ coefficients).
	If the fiber restriction $\HG X \lt \H X$ is surjective, 
	the action is called \defd{equivariantly formal}.	
	
	Note that $G$ itself admits a natural ``two-sided'' action of $G^2$
	given by $(g,h)\.x = gxh\-$.
	If $\defm{U} \leq G^2$ is any 
	closed, connected subgroup,
	%by the restriction of this action,
	we can equivalently write the homotopy quotient 
	of the restricted $U$-action in a more 
	symmetric way as
	\[
		G_U = \frac{EG \x G \x EG}{(eu,g,v\,\mn\mn f) 
		\sim (e,ugv,f)},
	\]
	where the relation ranges over $(u,v) \in U \leq G^2$ 
	and $(e,g,f) \in EG \x G \x EG$.
	In this case we call $G_U$ a \defd{homotopy biquotient}.
	In the case of greatest interest to us,
	where $U = \G \x H$ for some subgroups $\G,H \leq G$,
	we will also write 
	\[
	\defm{{}_\G G_H} \ceq EG\, \ox_\G \,G\, \ox_H \,EG.
	\]
	Since the right-multiplication action of $H$ on $G$ is free,
	the fiber of the map ${}_\G G_H \lt EG \ox_\G G/H = (G/H)_\G$
	given by $e \ox g \ox f \lmt e \ox gH$ is the contractible $EG$,
	so this map is a homotopy equivalence
	and $\H\mn({}_\G G_H) \iso \H_\G(G/H)$
	is the equivariant cohomology of the left $\G$-action
	by translation on 
	the homogeneous space $G/H$.
	The projection 
	$G_{\G \x H} \lt B(\G \x H)$
	of the Borel fibration
		can be equivalently written as 
	%\eqn{
	$	{}_\G G_H \lt B\G \x BH$
	%,\\
	%	e \ox g \ox f &\lmt \ (e\G, H\mn f)
	%}
	in this case, and the induced ring map as
%	inducing a ring map $\H_U \lt \H_U G$.
%	In case $U = \G \x H$,
%	the image of this map is $(e\G,H\mn f) \in B\G \x BH$
	$\H_\G \ox \H_H \lt \H_\G(G/H)$.
	In the event $\G = H$, 
	we economize on parentheses by writing
	$\defm{H_H^{\ox 2}} \ceq \H_H \ox \H_H$.
\end{definition}

\begin{construction}\label{rmk:Kapovitch}
	Given a two-sided action of $U \leq G^2$ on $G$,
	there is a map from the associated Borel fibration
	$G \to G_U \to BU$
	to the $G$-bundle
	\eqn{
		\defm{\D}\: BG \homeo EG \ox_G EG &\lt BG \x BG,\\
							e \ox f\ \ \ &\lmt\  (eG,G\mn f),
	}
	given on the total spaces by
%	
%	We will apply the theorem in the following situation.
%
%
%	and consider the homotopy quotient $\defm{G_U} \ceq E(G^2) \ox_U G$;
%	note $\KzGKz$, whose cohomology we aim to compute,
%	is $G_U$ for $U = K_0^2$.
%	The associated Borel fibration %can be modeled as
%	\eqn{
%		G_U = (EG)^2 \ox_U G	&\lt 	\xt U{EG}{EG} = BU,\\
%		(e,e') \ox g &\lmt 	e \ox e',
%	}
%	admits a $G$-bundle map to the bundle 
%	given on total spaces as 
	$(e,g,f)U \lmt eg \ox f = e \ox gf$.
	The induced map $BU \lt BG \x BG$
	downstairs is the map 
	$\defm{Bi}\: (e,f)U \lmt (eG,G\mn f)$
	induced by the inclusion 
	$\defm{i}\: U \longinc G$.
	\Cref{thm:FHTfibration}, applied to the bundle map
	\quation{\label{eq:pullback}
		\begin{gathered}
		\xymatrix@C=.75em@R=2.75em{
			G 	\ar@{=}[r]\ar[d]				&	G\ar[d]\\
			G_U \ar[r]\ar[d] 					& 	BG\ar[d]^\Delta\\
			BU \ar[r]_(.34){Bi}		&	BG \x BG,
		}
		\end{gathered}
	}
	will produce a model of $G_U$ 
	from models of the maps $Bi$ and $\Delta$.\footnote{\ 
		Kapovitch himself applied the result in the 
		case where $U$ acts freely on $G$,
		so that the homotopy biquotient $G_U$ 
		is homotopy equivalent to the naive orbit space $G/U$,
		which is in fact a smooth manifold, 
		conventionally called a \emph{biquotient}.
	}
	
	The first observation is that 
	by Hopf's theorem~\cite[Satz I]{hopf1941hopf},\footnote{\ 
		Hopf cares about the Pontrjagin ring $H_*(G;\R)$,
		the dual Hopf algebra to $\H(G;\R)$, 
		but the result is the same.
	}
	the cohomology of $G$ is a Hopf algebra, 
	whose underlying algebra is the exterior algebra $\eP$ on
	the space of primitive elements $\defm{P} = P\H G$.
	Thus $G$ admits a pure Sullivan model $(\ext P, 0)$ 
	with trivial differential.
	
	For the second observation, 
	write \defm{$\Sigma P$} for the \defd{suspension} of $P$,
	the graded vector space defined by $(\Sigma P)_n = P_{n-1}$,
	and $\defm Q = \defm{Q}\HG \ceq \wt H^*\mn(BG)/\wt H^*\mn(BG)^2$ for the so-called 
	\defd{module of irreducibles} 
	of $\HG$.
	This can be seen as a functorial version of the notion of a ``space of generators'': 
	a set of homogeneous $\Q$-algebra generators for $\HG$
	is exactly the image of a basis for $Q\HG$
	under a graded linear splitting of the projection $\HG \longepi Q\HG$.
	Borel's theorem~\cite[Thm.~19.1]{borelthesis} 
	states that the transgression
	in the \SSS of the universal bundle $G \to EG \to BG$
	induces a natural graded linear isomorphism 
	\[\defm\tau\:\Sigma P \isoto Q\HG,\]
	and that moreover, $\HG$ is the polynomial algebra on any 
	homogeneous lift of $Q\HG$.
	Picking such a lift at random, we can write $\tau P = \Sigma P = Q \HG \leq \HG$.
	The same of course holds for $\H_U$.
	
	Thus a model of $Bi\: BU \lt BG \x BG$ 
	is given by its own cohomology $\HG \ox \HG \lt \H_U$
	equipped with the zero differentials.
	The map $\D\: BG \lt BG \x BG$ can be identified
	up to homotopy with the diagonal,
	so it induces in cohomology the cup product $\HG \ox \HG\lt \HG$, 
	whose kernel is the ideal $(1 \ox \tau z - \tau z \ox 1: z \in P)$ 
	since $\HG$ is the symmetric algebra on $Q = Q\HG = \tau P$. 
	The transgression in the \SSS of the right bundle in (\ref{eq:pullback}) 
	thus sends each homogeneous $z \in P$ to $1 \ox \tau z - \tau z \ox 1$,
	suggesting we model $\D$ 
	by the inclusion \[(\HG \ox \HG,0) \longinc (\HG \ox \HG \ox \H G,d'),\]
	where $\defm{d'}$ vanishes on $\HG \ox \HG$
	and takes $z \in P$ to $1 \ox \tau z - \tau z \ox 1$%
	~{\cite[Prop., p.~157]{eschenburg1992biquotient}}.
	
	Then \Cref{thm:FHTfibration} yields the model
	\[
	(\H_U \ox \H G, d)
	\]
	for $G_U$,
	where $d$ is the derivation that vanishes on $\H_U$
	and takes $z \in P \H G$ to $(Bi)^*d'z \in \H_U$
	for $d'$ as in %the differential 
	%in the model %$(\HG \ox \HG \ox \H G,d')$ 
%	of %$BG$
	%in 
	the preceding paragraph.
	Specializing to $U = \G \x H$, 
	we get a pure Sullivan algebra
	\[
	(\H_{\G} \ox \H_H \ox \HG,d)
	\]
	modeling ${}_{\G}\mn G_{H}$, 
	where $d$ vanishes on $\H_\G \ox \H_{H}$ and 
	takes a primitive $z \in P \H G$
	to 
	\[
	(\rho_{\G}^* \ox \rho_{H}^*)(1 \ox \tau z - \tau z \ox 1) \ox 1,
	\]
	for $\defm{\rho_{\G}^*}\: \HG \lt \H_{\G}$
	and $\defm{\rho_{H}^*}\: \HG \lt \H_H$
	the maps induced 
	by the respective inclusions of $\G$ and $H$ in $G$.
	%Specializing again to the case $H = K_0$, 
	%we get a pure Sullivan algebra
	%\[
	%(\defm{\HKzt} \ox \HG,d) = (\HKz \ox \HKz \ox \H G, d)
	%\]
	%modeling $\KzGKz$, 
	%where $d$ vanishes on $\HKzt$ and 
	%takes a primitive $z \in P \H G$
	%to 
	%\[
	%(\rho^* \ox \rho^*)(1 \ox \tau z - \tau z \ox 1) \ox 1,
	%\]
	%for $\defm{\rho^*} = \rho^*_{K_0}$.
	Specializing to the case $H = 1$,
	we get back Borel's variant~\cite[\SS9]{cartan1950transgression}\cite[\SS25]{borelthesis}
	of the {Cartan algebra}
	computing $\H(G/K_0)$.
	%, which we will discuss more in 
	%\Cref{sec:culture}.
%	
%	An important algebraic observation about this model is that 
%	it can be considered as a Tor.
%	Namely, the transgression isomorphism $\tau\:P = P \H G \isoto Q\HG = Q$
%	induces a differential on $\HG \ox \H G \iso \ext Q \ox \ext P$
%	under which the cohomology is just $H^0 = \Q$;
%	this acyclic \CDGA is called the Koszul algebra
%	and can be seen as a model for the total space of the universal
%	bundle $G \to EG \to BG$.
%	The Koszul \CDGA
%	can be viewed as an $\HG$-module resolution $\HG \ox \ext P \to \Q$
%	of $\Q = \HG / H^{\geq 1}_G$,
%	where the resolution degree is given by 
%	$(\HG \ox \H G)^{-p} \ceq \HG \ox \ext^p P$.
%	By construction, the Kapovitch model $(\H_U \ox \H G,d)$
%	is isomorphic with the extension 
%	$\big(\H_U \ox_\HG (\HG \ox \H G), 1 \ox \tau\big)$
%	of the Koszul \CDGA,
%	so we conclude
%	\quation{\label{eq:KapTor}
%		\H(G_U) \iso \H(\H_U \ox \H G,d) \iso \Tor^*_{\HG}(\Q,\H_U).
%	}
\end{construction}

We will need a structure theorem on $\H(G/H)$
arising from an analysis of this model in the classical case $U = 1 \x H$.

\bthm[{\cite[Thm.~2, p.~211]{onishchik}}]\label{thm:formalCDGA}
Let $G$ be a compact, connected Lie group 
and $H$ a closed, connected subgroup. 
Then $G/H$ is formal if and only if there is an isomorphism
\[
\H(G/H) \iso (\H_H \ox_{\HG} \Q) \,\ox \ewP,
\]
where the exterior factor $\ewP$ is generated 
by a vector subspace $\wP$ of dimension $\rk G - \rk H$ 
and the kernel of $\dsp\H_H \lt \H_H \ox_{\HG} \Q$
can be generated by a regular sequence of elements of $\H_H$
lying in the image of 
$\smash{H^{\geq 1}_G \lt H^{\geq 1}_H}$.
%This in turn occurs
%if and only if $\H_H \ox_{\HG} \Q$ is a complete intersection ring.
\ethm

We can now prove \Cref{thm:cohomisotropyformal}.

\begin{proof}[Proof of \Cref{thm:cohomisotropyformal}]
We apply the Kapovitch model of \Cref{rmk:Kapovitch}
in the cases where $U$ is each of the groups 
\[	H \x H
\quad
>
\quad
1 \x H
\quad
>
\quad
1 \x 1.
\]	 
The corresponding homotopy quotients 
are the total spaces of the fibrations which form the columns of the diagram
%of the {\SSS}s of the fibrations (columns)
\[
\xymatrix{
	G \ar@{=}[r] \ar[d]	& G \ar@{=}[r]  \ar[d]	& G \ar[d]			\\
	{}_H G_{H}	\ar[d]	& G_H	\ar[l]	\ar[d]	& G \ar[l]\ar[d]	\\
	BH \x BH			& BH	\ar[l]			& {*}.\!\mn\mnn \ar[l]
}
\]
% $G \to G/H \simto G_{1 \x H} \to G_{H \x H} = EH \ox_H (G/H)$
The corresponding models form the sequence 
%G \to G_{H \x H} \to BH \x BH$
%and
%$G \to G_H \to BH$
%and 
%$G \to G \to *$.
%In total we model by
%maps
\[
(\H_H \ox \H_H \ox \H G,d) 
\xrightarrow{\e \,\mn\ox\,\mn \id \ox \id} 
(\H_H \ox \H G,\bar d)
\xrightarrow{\e \,\mn\ox\,\mn \id} 
(\H G, 0)
\]
of pure Sullivan algebras,
where $\e\: \H_H \to H^0_H = \Q$ is the augmentation and
$\defm d$ (resp., $\defm{\bar{d}}$) vanishes on $H_H^{\ox 2}$ (resp., $\H_H$) 
and takes a primitive element $z \in P \H G$
to 
\[
1 \ox \rho^*\tau z - \rho^*\tau z \ox 1 \in H_G^{\ox 2} \qquad 
(\mbox{resp., } \rho^*\tau z \in \HG).
\]
Here, as explained in the preceding discussion, $\rho = B(H \inc G)$
and $\tau\: P\H G \simto Q \HG \mono \HG$
is a lifting of the transgression in the \SSS of $G \to EG \to BG$.

	As $H$ acts equivariantly formally on $G/H$, we know
	\cite[Thm.~A]{carlsonfok2018}
	that $G/H$ is formal, which by \Cref{thm:formalCDGA} 
	means that 
	$\ewP \iso \im\mn\big(\H(G/H) \lt \H G \big)$
	is an exterior algebra on $(\rk G - \rk H)$ elements.
	By assumption, $\H_H(G/H) \lt \H(G/H)$ is surjective,
	so the image of $\H_H(G/H) \lt \H G$ is also $\ewP$.
	Writing $\defm{\cP}$ for a complement of $\wP$ in $P\H G$,
	and picking cocycles %in the first algebra 
	representing elements that map to $\ewP$,
	we may factor $(H_H^{\ox 2} \ox \H G,d)$ as
	$(H_H^{\ox 2} \ox \ecP,d) \ox {(\ewP,0)}$.
	On applying $\e \ox \id^{\ox 2}$,
	we recover the factorization
	$(\H_H \ox \ecP,\bar d) \ox {(\ewP,0)}$
	of $(H_H^* \ox \H G,\bar d)$ 
	figuring in the proof of \Cref{thm:formalCDGA}.
	Since $G/H$ is formal,
	%again by
	 \Cref{thm:formalCDGA} shows
%	\[
%	\H(\H_H \ox \ecP) \iso \im\mn\big(\H_H \to \H(G/H)\big) \iso \H_H /(\bar d\cP)\]
%	is a complete intersection ring,
%	or in other words, 
	a homogeneous basis $(\defm{x_j})$ of $\smash{\bar d\cP}$ forms a regular sequence in $\H_H$.
	Particularly, $\H_H$ is free over 
	the polynomial subring $\defm C \ceq \ext[\bar d\cP] = \Q[x_j]$.
	The nonzero part of $d$ is determined its restriction
	to ${\cP}$,
	which factors as
	\[
	\cP \isoto \bar d\cP \longinc C \os{\defm \nu}\lt C^{\ox 2} \longinc H_H^{\ox 2},
	\]
	where $\nu$ is given
	%by 
	%Particularly, $\H_H$ is a finite free module over $\Q[d\cP]$~\cite[XXX]{neuselsmith}.
	%The map
	%$\H_H \longmono \H_H \ox \H_H$
	%given on a set of polynomial generators $x$ 
	by $x_j \lmt 1 \ox x_j- x_j\ox 1$.
	Now $\nu$ makes $C^{\ox 2} = \Q[x_j]^{\ox 2}$ 
	a free module (indeed, a polynomial algebra) over $C = \Q[x_j]$, 
	since
	%makes the former free over the latter since 
	$1 \ox x_j- x_j\ox 1$ and $1 \ox x_j+ x_j\ox 1$
	are algebra generators for the former,
	so $H_H^{\ox 2}$ is free over $C$ as well.
	But then since $(C \ox \ecP,\bar d)
	= (\ext [\bar d\cP] \ox \ext \cP,\bar d)$
	is acyclic,
		\[
	\H(H_H^{\ox 2} \,\ox \ecP,d) =
	H_H^{\ox 2} \,\ox_C \,\H(C \ox \ecP,\bar d) \,\iso\, 
	H_H^{\ox 2} \,\ox_C \,\Q \,\iso\,
	%\frac{\H_H \ox \H_H}{(1 \ox x - x \ox 1 : x \in d\cP)}
	\xt{C}{\H_H\,}{\,\H_H} \,\iso\, 
	\xt{\HG}{\H_H}{\H_H}.\qedhere
		\]
	
We then have more specific structure theorems as corollaries
in the special cases mentioned in the introduction.

\bcor[{\cite[Prop.~68, p.~161]{kumarvergne1993equivariant}\cite[Thm.~11]{tu2010homogeneous}}]\label{thm:equalrank}
Let $G$ be a connected Lie group, 
$H$ a closed subgroup of equal rank,
and $\G$ another closed subgroup.
Then there is a ring isomorphism
\[
\H_\G(G/H;\Q) \iso \xt{\H(BG;\Q)}{\H(B\G;\Q)}{\H(BH;\Q)}.
\]
\ecor
\bpf
If $\G = H = T$ is a maximal torus of $G$, 
this follows from \Cref{thm:cohomisotropyformal}
since $\rk G - \rk T = 0$.
If $\G$ and $H$ are connected subgroups containing $T$,
then 
\[
	\H_{\G \x H} G 
		\iso 
	\H_{T \x T}(G)^{W_\G \x W_H}
		\iso
	(\HT \ox_\HG \HT)^{W_\G \x W_H}
		\iso
	\H_\G \ox_\HG \H_H
\]
by a well-known lemma~\cite[Prop.~III.1, p.~38]{hsiang}, 
where $\defm W$ denotes the Weyl group.
To pass to disconnected subgroups $\G,H$
with respective identity components $\G_0,H_0$,
note that ${}_{\G_0 } \mn G_{H_0} \lt {}_{\G} G_{H}$
is a principal ($\pi_0 \G \x \pi_0 H$)-bundle,
so a standard lemma on invariants~\cite[Prop.~3G.1]{hatcher} 
lets us identify
\[
\H({}_{\G} G_{H}) \iso
\H({}_{\G_0} \mn G_{H_0})^{\pi_0 \G \x \pi_0 H} \iso
(\H_{\G_0} \ox_\HG \H_{H_0})^{\pi_0 \G \x \pi_0 H} \iso
\H_\G \ox_\HG  \H_H.
\]

To pass to the general case, note that 
by the above in the case $G = \G$, 
we have 
\[ 
\HG(G/H) =\HG \ox_\HG \H_H \isoto \H_H \longepi \Q \ox_\HG \H_H \isoto \H(G/H)
\]
surjective. 
For any other $\G < G$, we can factor 
$G/H \lt (G/H)_G$
as $EG \x G/H \to EG \ox_\G G/H \to EG \ox_G G/H$,
so $\H_\G(G/H) \lt \H(G/H)$ is surjective as well.
Thus the \SSS of $G/H \to (G/H)_\G \to B\G$
collapses at the $E_2$ page,
so that on the level of $\H_\G$-modules one has
$\H_\G(G/H) \iso \H_\G \ox \H(G/H)$.
This implies in turn that
the map $\H_\G \ox \H_H \lt \H({}_\G G_H)$
induced by ${}_\G G_H \lt B\G \x BH$ is surjective.
From the \SSS map induced by \Cref{eq:pullback} 
it follows the image can be identified with $\H_\G \ox_\HG \H_H$
as claimed.
\epf
	
\bcor\label{thm:Hsurj}
Let $G$ be a connected Lie group and $H$ a closed, 
connected subgroup such that 
$\H G \lt \H H$ is surjective. Then
\[
\H_H(G/H) \iso \H_H \ox \ewP,
\]
where $\ewP$ is isomorphic to the image of $\H(G/H;\Q) \lt \H (G;\Q)$.
Particularly, $\dim_\Q \wP = \rk G - \rk H$. 
\ecor
\bpf
In this case, 
a theorem of Borel\footnote{\ 
	seen for example by taking $EH = EG$ 
and following the
transgressions in the map of spectral
sequences induced by the bundle map
$(H \to EH \to BH) \lt (G \to EG \to BG)$}
shows $\HG \lt \H_H$ is surjective as well,
so that $\dsp\smash{\H_H \ox_\HG \H_H} = \H_H \ox_{\H_H} \H_H \iso \H_H$.
\epf

\bcor\label{thm:circle}
Let $G$ be a connected Lie group containing a subgroup $H \iso S^1$
such $N_G(H)$ has two components.
Then
\[
	\phantom{,\qquad |\s| = |s| = 2, \quad |z| = 3 }
	\H_H(G/H) 
		\iso 
	\frac{\Q[\s,s]}{(\s^2 - s^2)} \ox \frac {\H G}{(z)},
	\qquad |\s| = |s| = 2, \quad |z| = 3,
\]
where $z$ is represented in de Rham cohomology by the left-invariant 
\emph{Cartan 3-form} given at 
$1 \in G$ by $(u,v,w) \lmt \big\langle u,[v,w]\big\rangle$
for $\langle-\rangle$ an $\Ad(G)$-invariant inner product on 
the Lie algebra $\fg$ of $G$.
\ecor
\bpf
	In this case, $N_G(H)$ contains an element $g$
	such that $ghg\- = h\-$ for $h \in H$
	as this is the only nontrivial continuous automorphism of $S^1$.
	If we choose a generator $\s$ of $\H_H = \H_{S^1} = \Q[\s]$,
	then the image of $\HG \lt \H_H$ must lie in 
	$(\H_H)^{N_G(H)} = \Q[\s^2]$.
	The proof of equivariant formality in this case~\cite[Prop.~1.4]{carlson2018eqftorus}
	shows $\s^2$ is indeed in the image of $\HG$
	and is the transgression of $z \in H^3 G$
	in the \SSS of $G \to G/H \to BH$,
	so $\smash{\H_H \ox_\HG \H_H \iso 
		\Q[\s] \ox_{\Q[\s^2]} \Q[\s]}$
	is the first factor in the statement of \Cref{thm:cohomisotropyformal}
	and $\wP$ is spanned 
	by a complement of $\Q z$ in $P\H G$.
\epf

\begin{EG}
	%
	%\section{Other examples}
	
	To illustrate the versatility
	of the model of \Cref{rmk:Kapovitch},
	%lest the reader imagine it be limited to real Grassmannians,
	we briefly list some representative examples.
	All the following are computable with equal facility
	using the model $(\H_{\G} \ox \H_H \ox \H G,d)$
	once one has expressions for the maps $\rho_\G^*$ and $\rho_H^*$.
	The spaces in all of these examples but the last two
	are formal,
	and all but the last three 
	can be derived from \Cref{thm:cohomisotropyformal}
	or its corollaries.

	\eqn{
		%\H_{\U(k) \x \U(n)} \Z\Big(\quotientmed{\U(n+k)}{\U(k) \x \U(n)};\Z\Big)
		%				&= 	 \Z[c,c',\vk,\vk']/(cc'-\vk\vk');
		%				\\
		\H_{\Sp(k) \x \Sp(n)} \Big(\quotientmed{\Sp(n+k)}{\Sp(k) \x \Sp(n);\Z}\Big)
		&= 	 \Z[p,p',\pi,\pi']/(pp'-\pi\pi');
		\\
		\H_{\U(n)}\big(\Sp(n)/\U(n);\Z\big)
		&=	 \Z[c,\vk]/(c\bar c - \vk \bar \vk);
		\\
		\H_{G_2}\big(\Spin(9)/G_2\big) 
		&= \Q[y_4,y_{12}] \ox \ext[z_7,z_{15}];
		\\
		\H_{G_2}(F_4/G_2) &= \Q[y_4,y_{12}]\ox\ext [z_{15},z_{23}];
		\\
		%\H_{F_4}(E_6/F_4) &= \Q[y_4,y_{12},y_{16},y_{24}]\ox\ext[z_{9},z_{17}];
		%					\\
%		\H_{E_6}(E_7/E_6) &= \Q[y_4,y_{10},y'_{10},y_{12},y_{16},y_{18},y'_{18},y_{24},
%		]\ox\ext[z_{19},z_{27},z_{35}];
%		\\
		\H_{\U(2k)}\big(\U(2n)/\Sp(n);\Z\big)
		&=	 \frac{\Z[c_2,\ldots,c_{2k},p_1,\ldots,p_n]}{
			(c_{2j} - p_j)} 
		\ox \ext[z_{2k+1},z_{2k+3},\ldots,z_{2n-1}];
		\\
		\H_{\SU(3) \x \SU(3)}\Big(\quotientmed{\SU(6)}{\SU(3) \x \SU(3)}\Big)
		&=	 \Q[c_2,c_3,c'_2,c'_3,
		\vk_2,\vk_3,\vk'_2,\vk'_3] 
		/ (cc' - \vk\vk');	\\
		\H_{S^1}\big(\SU(3)/S^1\big) 
		&= 	\, 	\frac{\Q[v_2,w_2] \ox \ext[y_7]}
		{(vw,yw,w^3)}, \qquad 
		S^1 = \big\{\!\diag(z,z,z^{-2})\big\}.
	}
	
	\medskip
	
	\nd Here $c = 1 + c_1 + \cdots + c_n,c',\vk,\vk'$ are total Chern classes
	and 
	$\bar c = 1 - c_1 + \cdots + (-1)^n c_n,\bar\vk$ their duals,
	$p$, etc., are symplectic Pontrjagin classes,
	$\Z[p] \ceq \Z[p_1,p_2,\ldots]$ and so on,
	and elements that are not Chern or symplectic Pontrjagin classes
	are indexed by their degrees.

\end{EG}

\end{proof}

\begin{history}\label{rmk:Borel}
There is a construction of the Kapovitch model using less explicit
rational homotopy theory, 
which the author found by generalizing
the construction of the Cartan model in Borel's thesis
before he realized this model was already known.
(It in fact
already follows from Borel's thesis~\cite[Thm.~24.1', p.~185]{borelthesis}
that this will work.)
The idea of this proof is to note that \eqref{eq:pullback}
induces a map of Serre spectral sequences,
and that the transgressions in the right sequence,
$z \lmt 1 \ox \tau z - \tau z \ox 1$ for $z \in P \H G$,
are known, hence determining the transgressions in the left 
sequence. 
One needs a functorial \CDGA model $M(-)$,\footnote{\ 
Borel uses a real model based in sheaf theory and ultimately 
induced from the de Rham algebra; 
the author used the algebra $\APL$ of polynomial differential forms.}
and then lifts generators $z$ of $\H G$ to cochains $\zeta$ in $M(G_U)$
which on restriction to $M(G)$ 
represent generators of $\H G$.
As the classes of these $\z$ transgress in the \SSS,
their coboundaries $d_{M(G_U)} \z$ lie in the image of $M(BU) \lt M(G_U)$,
and picking preimages $\d z$ in $M(BU)$,
one can uniquely define a derivation $\d$ on $M(BU) \ox \H(G)$
taking $z \lmt \d z$ and restricting to $d_{M(BU)}$ on the $M(BU)$
factor.
There is then a natural \CDGA map $\big(M(BU) \ox \H G,\d\big) \lt M(G_U)$;
the lifting $z \lmt \z$
extends to a homomorphism because $M$ is commutative.
The map of spectral sequences induced by the 
induced by the skeletal
filtration of $BU$
is easily seen to be an isomorphism on $E_2$ pages.
It follows the \CDGA map is a quasi-isomorphism.
Since $M$ is commutative
and $H_U$ is a polynomial ring,
one can lift it to a polynomial subring of $M(BU)$
simply by lifting generators and extending uniquely,
which determines a map of 
the Kapovitch model $\H_U \ox \H G$
into $M(BU) \ox \H G$. 
Again there is a map of filtration spectral 
sequences showing that this is a quasi-isomorphism
and thus that the Kapovitch model indeed computes $\H(G_U)$.
%The skeletal filtration of $BU$ induces the filtration of 
%the Kapovitch model and the standard filtration of $G_U$
%inducing the \SSS,
%so our \CDGA map induces a map from the filtration spectral
%sequence to the Serre spectral sequence of $G \to G_U \to BU$,
%which it is easy to see is already an isomorphism on $E_2$ pages.
%It follows our \CDGA map is a quasi-isomorphism,
%so the cohomology of the Kapovitch model is indeed $\H(G_U)$.
The case $U = 1 \x H$ is Borel's argument that the Cartan model
computes $\H(G/H)$.

There is a second way to the model. 
Borel's theorem can be seen as showing that $EG$ 
is modeled by $(\HG \ox \H G,\tau)$,
where $\tau$ is the derivation
defined as a lifting of the transgression on $P \H G$
and vanishing on $\H G$.
In other words, this \CDGA, often called a \emph{Koszul complex},
is acyclic. Grading it by exterior degree, 
one can view it as an $\HG$-module resolution of $\Q$
and so use it to define Tor.
%Then the cohomology of the Kapovitch model can 
%be identified with $\Tor^*_{\HG}(\Q,\H_U)$. 
Similarly, the model $(H_G^{\ox 2} \ox \H G,d')$
from \Cref{rmk:Kapovitch} can be viewed as an 
$H_G^{\ox 2}$-module resolution of $\HG$.
Then one can compute $\Tor_{H_G^{\ox 2}}^*(\H_G,\H_U)$
as the cohomology of 
$(H_G^{\ox 2} \ox \H G) \ox_{H_G^{\ox 2}} \H_U
= \H_U \ox \H G$, which is the Kapovitch model.
On the other hand, this Tor is the $E_2$
page of the 
\EMSS associated to the homotopy pullback square \eqref{eq:pullback},
which therefore collapses.
This approach, however, does not show 
this Tor has the correct ring structure.

In the case $U = 1 \x H$, 
this approach to computing $\H(G/H;k)$
over a suitable coefficient ring $k$
is classical and first 
appeared in Paul Baum's thesis
as the claim that 
\quation{\label{eq:Baum}
	\H(G/H;k) \iso \Tor^*_{\H(BG;k)}\big(k,\H(BH;k)\big)
}
when $k$ is a field 
and $\H(G;k)$ and $\H(H;k)$ are exterior algebras on odd-degree generators.
A mistake however resulted in the claim only being 
only being demonstrated when $\Tor^{\geq 3}$ vanishes
in the published version~\cite{baum1968homogeneous}.
All work aimed at recovering the claimed result
involves higher homotopy structures such as $A_\infty$-algebras%
~\cite{gugenheimmay,munkholm1974emss,wolf1977homogeneous}.
The best such result is arguably that \eqref{eq:Baum} holds
when $k$ is a principal ideal domain
such that $\H(BG;k)$ and $\H(BH;k)$ are polynomial algebras 
and the Steenrod square $\Sq^1$ vanishes on $\H(BH;k)$ 
in case $k$ has characteristic two.
The ring structure, except in the case $G$ and $H$ have equal rank~\cite[Cor.~7.5]{baum1968homogeneous},
was unaddressed until a series of
papers by Matthias Franz~\cite{franz2019homogeneous},
released since the present work was submitted,
showed that \eqref{eq:Baum} 
is in fact an isomorphism of graded $k$-algebras
for any principal ideal domain $k$ for which $2$ is invertible and 
$\H(BG;k)$ and $\H(BH;k)$ are polynomial rings.

The \EMSS approach can be generalized in a different direction, 
in the case $U = \G \x H$,
by instead considering the pullback square
\[
\xymatrix@C=1em{
	{}_\G G_H\ar[r]\ar[d]&BH\ar[d] \\
	B\G	\ar[r]		&\,BG	.
}
\]
Singhof~\cite{singhof1993} showed using the above-cited work of Munkholm~\cite{munkholm1974emss}
that if the ($\G \x H$)-action on $G$ is free,
meaning ${}_\G G_H$ is homotopy equivalent to the honest biquotient
$\G \lq G /H$,
then the associated \EMSS collapses at 
$E_2 = \Tor^*_{\H(BG;k)}\big(\H(B\G;k),\H(BH;k)\big)$
for any coefficient ring $k$ 
such that the rings $\H(B\G;k)$, $\H(BH;k)$, and $\H(BG;k)$
are all polynomial on even-degree generators,
though again this only partially determines the ring structure
of $\H(\G \lq G /H;k)$.
In the light of Franz's recent work it now seems it may be possible 
to extend these and the present results 
to prove the same ring structure for more general coefficients.

%\bprop
%In the case $G = \SO(2n+2k+\a+\b)$ for $0 \leq \a \leq \b \leq 1$
%and $K_0 = \SO(2k+\a) \x \SO(2n+\b)$,
%one may actually take $\Zt$ coefficients in the Kapovitch model of
%\Cref{rmk:Kapovitch}.
%\eprop
%\bpf
%In the diagram \eqref{eq:pullback} 
%with $U = K_0 \x K_0$, 
%the cohomology of each of $G$, $BG$, $BG$, and $BK_0$
%is a free module over $\Zt$.
%so it follows from the \SSS of the left bundle
%that 
%\epf
%
%The Kapovitch model //
%for $U = 1 \x H$ (i.e., for $G/H$)
%is due (in a rather different derivation with $\R$ coefficients) 
%to Cartan~\cite[Thm.~5, p.~216]{cartan1950transgression}, 
%building on work of Weil and Chevalley.
%%
%%the comparison with the \SSS of $G \to G/H \to BH$,
%%which is the ``pullback argument'' from the universal bundle
%%mentioned in the preceding footnote,
%%is essentially due to Borel in his thesis%
%%~\cite[Thm.~25.2]{borelthesis}.
%The isomorphism \eqref{eq:KapTor}
%in this case, namely
%$\smash{\H(G/H) \iso \Tor^*_{\HG}(\Q,\H_H)}$,
%means Cartan's theorem that $(\H_H \ox \H G,d)$
%is a model can viewed as a strong
%collapse theorem for the \EMSS of the fibration $G/K \to BK \to BG$.
%A proof of Cartan's theorem by demonstrating this \EMSS
%collapses under certain conditions 
%is the central result of the thesis of Paul Baum~\cite{baumthesis}.
\end{history}

\section{Grassmannians}\label{sec:Grass}

We come to our motivating example.
%\begin{notation}
%We will need more notation to state the results efficiently.

%to propagandize these techniques.
%in the hope they might become better known.

%The two goals of this note, to be short and to be expository,
%are in tension, much to the expense of the former. 
%While we have been able to make the proof itself quite brief,
%to explain 
%why the model that we find so convenient should exist at all 
%takes about two pages.
%%We think in light of its great utility, 
%%this loss is acceptable.
%We round the note out with an alternate proof and 
%some historical remarks.
%The only potentially original contribution is
%the observation that the model we use 
%in fact applies.

%\medskip

\begin{notation}\label{rmk:bignotation}
Let $\defm{G_\ell(\R^{\ell+m})}$
	denote the Grassmann manifold of unoriented $\ell$-dimensional real subspaces
	of $\R^{\ell+m}$
	and $\defm{\wt G_{\ell}(\R^{\ell+m})}$
	that of oriented subspaces.
	%\end{notation}
We write $\defm{\O(\ell)}$ for the orthogonal group and 
$\defm{\SO(\ell)}$
for the special orthogonal group.
%Given a group $G$,
%we adopt the abbreviations $\defm{\HG} \ceq \defm \HG(\pt) = \H BG$
%for the cohomology of the classifying space
%and $\defm{H_G^{\ox 2}} \ceq \HG \ox \HG$.
The total Pontrjagin class in $\H_{\O(\ell)}$ 
is the sum $1 + p_1 + \cdots + \smash{p_{\floor{\ell/2}}}$ 
of the Pontrjagin classes $\smash{{p_j} \in {H^{4j}_{\O(\ell)}}}$
and the ring $\H_{\O(\ell)}$ is the symmetric algebra on the $p_j$.
When $\ell$ is odd, the map $B\SO(\ell) \lt B\O(\ell)$ induces an isomorphism $\smash{\H_{\O(\ell)} \isoto \H_{\SO(\ell)}}$.\footnote{\ 
	For both these statements, 
	it is critical that $2$ is invertible in $\Q$;
	there are many more torsion classes of order $2$
	in $\H\big(B\O(\ell);\Z\big)$ and $\H\big(B\SO(\ell);\Z\big)$ 
	which arise from
	Bockstein images of Stiefel--Whitney classes 
	in the integral cohomology ring.	
}
When $\ell$ is even, the top Pontrjagin class $p_{\ell/2}$  
is also the square $\smash{e^2}$ of the Euler class
${e} \in \smash{H^{\ell}_{\SO(\ell)}}$---%
we adopt the convention that $e$ is defined but zero when $\ell$ is odd---%
%Definitionally, $p_0 = 1$.
and~\cite[Thm.~3.16]{VBKT}
\[\H_{\SO(\ell)} \iso \Z[\sfrac 1{\,2}][p_1,\ldots,p_{(\ell/2)-1},e].\]

Let $0 \leq \defm{\a} \leq \defm{\b} \leq 1$ 
and $\defm k,\defm n \geq 0$ be natural numbers, and set 
\eqn{
	\defm{G} &\ceq \SO(2k+2n+\a+\b);\\
	\defm{K} &\ceq \mathrm S \big(\O(2k+\a) \x \O(2n+\b)\big);\\
	\defm{K_0} &\ceq \SO(2k + \a) \x \SO(2n+\b);\\
	\defm{T} &\ceq \SO(2)^{\oplus k} \x \SO(2)^{\oplus n},
}
the last being the ``diagonal'' maximal torus of $K_0$.
Then the right quotient $G/K = G_{2k+\a}(\R^{2k+2n+\a+\b})$
is the unoriented real Grassmannian
and $G/K_0 = \wt G_{2k+\a}(\R^{2k+2n+\a+\b})$
the oriented Grassmannian.
Only when $\a = \b = 1$ is the Grassmannian odd-dimensional;
we will call this the \defd{odd case} 
and the others the \defd{even cases}.
We aim to compute the rings $\H_\G(G/H)$ 
%for $\G,H < G$
for $\G,H \in \{K,K_0\}$;
the case $\G = T$ that He discusses follows from the case $\G = K_0$
by a standard lemma, as elaborated in \Cref{rmk:recovery}.

The natural maps $\H_\G \ox \H_H \lt \H_\G(G/H)$ 
discussed in \Cref{def:actions}
will be the source of most of the
classes in the statement of the theorem.
%makes $\KzGKz$ a $G$-bundle over $BK_0 \x BK_0$.
For $(\G,H) = (K_0,K_0)$, write respectively
$\defm{\pi,\pi',p,p'}$ for the total Pontrjagin classes
in
\[
\H_{\O(2k + \a)} < \H_\G,\qquad
\H_{\O(2n + \b)} < \H_\G,\qquad
\H_{\O(2k + \a)} < \H_H,\qquad
\H_{\O(2n + \b)} < \H_H.
\] 
Similarly, if $\G$ or $H$ is $K_0$, 
write respectively
$\defm{\e,\e',e,e'}$ for the Euler classes (zero if $\a$, resp. $\b$, is $1$)
of
\[
H^{2k}_{\SO(2k + \a)} < \H_\G,\qquad
H^{2n}_{\SO(2n + \b)} < \H_\G,\qquad
H^{2k}_{\SO(2k + \a)} < \H_H,\qquad
H^{2n}_{\SO(2k + \b)} < \H_H.
\]
We use the same notations for the images of these classes in 
$\HK \ox \HK$, in $\HKz \ox \HK$, and in $\HKz \ox \HKz$,
and for the further
images in $\HK(G/K)$, in $\HKz(G/K)$, and in $\HKz(G/K_0)$
under the natural maps.

Faced with an inhomogeneous graded ring element $x$ 
with homogeneous components $x_n$
(e.g., $p$),
we make the nonstandard abbreviation $\defm{\Zt[x]} \ceq \Zt[x_1,x_2,\ldots]$,
and similarly we denote by $\defm{(x)} \ceq (x_n)_{n \geq 0}$ 
the ideal generated by the {homogeneous components} of $x$.
We write $\defm{\ext [\eta]} $ for the exterior algebra 
%(\defd{\fontfamily{qpl}\selectfont \textsc{cga}})
%on a positively-graded rational vector space
over $\Z$
generated by
%and $\defm{\ext[\eta]}$ if $P$ is spanned by 
one element $\h$.
Specifically, in the odd case, we will be interested in 
the element $\defm \eta \in \smash{H^{2k+2n+1}_{K_0}(G/K_0)}$ which
restricts under $G \lt \KzGKz$ 
to that primitive generator of $\H G$
which transgresses in the \SSS of $G \to EG \to BG$
to the Euler class in $H^{2k+2n+2}_{\SO(2k+2n+2)}$.
Finally, we set
\[
\defm{\pont} \ceq \frac{\Zt[\pi,\pi',p,p']}{(\pi\pi'-pp')}.
\]
\end{notation}

\bthm[He, 2016]\label{thm:He}
The Borel equivariant cohomology rings
of the left actions of $K_0$ on
the oriented real Grassmannian
and of $K_0$ and $K$
%$\SO(2k+\a) \x \SO(2n+\b)$ and 
%$\mathrm{S}(\O(2k+\a) \x \O(2n+\b))$
on the unoriented Grassmannian 
%$G/K_0$ 
%and of $K$ and $K_0$ on the unoriented Grassmannian $G/K$
are as follows.

\vspace{-2ex}

\[\arraycolsep=12pt\def\arraystretch{2.5}
\begin{array}{r|lll}
(\a,\b)
& (1,1)\quad 
&(0,1) \quad
&(0,0)
\\
\hline
\H_{\SO(2k+\a) \x \SO(2n+\b)} \big(\wt G_{2k+\a}(\R^{2k+2n+\a+\b});\Zt\big) 
&\pont
\ox \ext[\h]\phantom{X}
&\pont[\e,e]\phantom{X}
&\dsp\frac{\pont[\e,\e',e,e']}{(\e\e' - ee')}
\\
\H_{\SO(2k+\a) \x \SO(2n+\b)} \big(G_{2k+\a}(\R^{2k+2n+\a+\b}) ;\Zt\big)
&\pont
\ox \ext[\h]
&\pont[\e]
&\dsp\frac{\pont[\e,\e',ee']}
{(\e\e' - ee')}
\\
\H_{\mathrm{S}(\O(2k+\a) \x \O(2n+\b))} \big(G_{2k+\a}(\R^{2k+2n+\a+\b}) ;\Zt\big)
&\pont
\ox \ext[\h]
&\pont
&\dsp\frac{\pont[\e\e',ee']}
{(\e\e' - ee')}
\end{array}
\]
%\tc{Check}

\nd 
Here the relations 
	$e^2 = p_k$, \ \
	$(e')^2 = p'_n$,\ \ 
	$\e^2 = \pi_k$, \ \  
	$(\e')^2 = \pi'_n$
are tacit
and the natural maps induced by coverings preserve
elements with identical names.\footnote{\ 
	It is possible  all the torsion is in fact of order $2$,
	and if so the statement will remain true if one replaces 
	the cohomology with integral cohomology modulo torsion of order $2$ 
	on the left-hand side, 
	and $\Zt$ with $\Z$ %in the ring expressions 
	on the right.}
%where we view
%$\HK$
%as the subalgebra of $\Z/2$-invariants
%of $\HKz$ using the double covering $BK_0 \lt BK$.
%and as before $(ee')^2 = p_n p_k$ and $(\e\e')^2 = \pi_n\pi_k$
\ethm

We recover the singular cohomology rings as a corollary.
Write $\dsp\defm{\redpont} \ceq 
\frac{\Zt[p,p']}{(pp'-1)}$.

\medskip

%It remains to obtain our statements about coefficients.
%
%\bpf[Proof of \Cref{thm:coefficients}]
%For any space $X$, we have a sequence of coefficient maps
%\[
%	\H(X;\Z) \lt 
%	\H(X;\Z) / 2\mbox{-torsion} \lt
%	\H\big(X;)\Z[\tfrac 1 2]\big) \lt
%	\H(X;\Q).
%\]
%We already know the Pontrjagin and Euler classes 
%in \Cref{thm:He}
%are integral,
%so our proposition amounts to the claim
%the only torsion in $\KzGKz$ and $\KGK$ is $2$-torsion.
%To that end, we look at the map of {\SSS}s
%we used to set up our model.
%Recall that
%all torsion of $\SO(n)$, of $B\O(n)$, and of $B\SO(n)$
%is $2$-torsion.
%It follows that there is only $2$-torsion
%in the $E_2$ page of the \SSS
%of $G \to \KzGKz \to BK_0 \x BK_0$.
%The map of Serre spectral sequences 
%is a chain map on each page.
%Since the kernel and cokernel
%in each transgression in the \SSS of $G \to EG \to BG$ is zero
%
%\epf

\bcor[{Leray~\cite{lerayCR1949e}, Borel~\cite[p.~192]{borelthesis}, 
	Takeuchi~\cite[p.~320]{takeuchi1962pontrjagin}}]\label{thm:cohomhomog}
The singular cohomology rings
of the oriented and unoriented real Grassmannians
are as follows.

%\eqn{
%	\H
%	\wt G_{2k}(\R^{2k+2n}) 
%		&\iso
%	\frac{\Q[p,p',e,e']}
%			{(pp' - 1,ee')},\\
%	\H
%		\wt G_{2k}(\R^{2k+2n+1}) 
%			&\iso 
%	\frac{\Q[p,p',e]}
%			{(pp'-1)},\\
%	\H
%	G_{2k+1}(\R^{2k+2n+2}) 
%			\iso 	
%	\H
%	\wt G_{2k+1}(\R^{2k+2n+2}) 
%			&\iso 
%	\frac{\Q[p,p']}
%				{(pp'-1)}
%			\,\ox\,
%		\ext[\eta],\\
%	\H
%	G_{2k}(\R^{2k+2n}) 
%			\iso				
%			%	\\	
%	\H
%	G_{2k}(\R^{2k+2n+1}) 									
%			&\iso 
%			\frac{\Q[p,p']}
%					{(pp' - 1)}.
%}
%
%
\vspace{-2ex}

\[\arraycolsep=12pt\def\arraystretch{2.5}
\begin{array}{r|lll}
(\a,\b)
& (1,1)\quad 
&(0,1) \quad
&(0,0)
\\
\hline
\H\big(\wt G_{2k+\a}(\R^{2k+2n+\a+\b});\Zt\big)
&\redpont
\ox \ext[\h]\phantom{X}
&\redpont[e]\phantom{X}
&\dsp\frac{\redpont[e,e']}
{(ee')}
\\
\H\big(G_{2k+\a}(\R^{2k+2n+\a+\b});\Zt\big)
&\redpont
\ox \ext[\h]
&\redpont
&\redpont
\end{array}
\]

\nd 
Here 
%$\eta$ restricts under $G \lt G/K_0$ 
%to the suspension in $H^{2k+2n+1} \SO(2k+2n+2)$ 
%of the Euler class in $H^{2k+2n+2}_{\SO(2k+2n+2)}$
%and 
the relations 
$e^2 = p_k$ and $(e')^2 = p'_n$
are tacit
and the natural map induced by 
a double covering preserves
elements with identical names.
\ecor
\bpf
Since the actions in question are
equivariantly formal,
in the bundles
$G/K_0 \to (G/K_0)_{K_0} \to BK_0$
and $G/K \to  (G/K)_{K_0} \to BK_0$,
the fiber restrictions are surjections
whose kernel is the ideal 
generated by the classes of positive degree
in the image of $\H(BK_0) \lt \HKz(G/K_0)$
and $\H(BK_0) \lt \HKz(G/K)$ respectively.
Thus the expressions above
follow from \Cref{thm:He}
on quotienting out those generators 
amongst $\pi - 1, \pi'-1, \e,\e'$, and $\e\e'$
which are relevant. 
\epf

\brmk\label{rmk:recovery}
The equivariant cohomology in He's original statement of the theorem
is recovered via the 
standard natural isomorphism~{\cite[Prop.~III.1, p.~38]{hsiang}}
\[
	\smash{\HT(X) \iso \xt {\H_{K_0}}\HT{\H_{K_0}(X)}}
\]
once we have the explicit descriptions of the injections $\HK \inc \HKz \inc \HT$ from the following discussion.
\ermk
\begin{discussion}\label{rmk:restrictions}
%To understand $\HK \longinc \HKz$,
%first note that $K_0$ is of index two in $K$, it is particularly
%normal, so $K/K_0$ is a group isomorphic to $\Z/2$.
%We may view $BK_0$ as the quotient $EK/K_0$,
%which admits a right action of the group $K/K_0$.
%With these identifications, 
%$BK_0 \lt BK$ is then a principal $\Z/2$-bundle covering.
%As $1/2$ is invertible in $\Zt$,
%by a standard lemma on invariants~\cite[Prop.~3G.1]{hatcher}
%we may identify $\HK$ with the $\Z/2$-invariants of $\HKz$.
To understand the maps $\HK \to \HKz \to \HT$,
we make the general observation that if  
$\G$ is a compact Lie group with maximal torus $T$, 
not necessarily connected,
we may identify $\H_\G$ 
with the subring of invariants of $\HT$
under $N_\G(T)$.
To see this, first note the identity component $\G_0$ is normal in $\G$,
so given any $\g \in \G$ one has $\g T \g\-$ a maximal torus of $\G_0$,
hence there is $\g_0 \in \G_0$
with $\g_0\g T \g\- \g_0\- = T$,
and each component of $\G$ contains an element of $N_\G(T)$.
There follows the diffeomorphism $\G/N_\G(T) \homeo \G_0/N_{\G_0}(T)$
of homogeneous spaces.
The later space is known to be rationally acyclic~\cite[Lem.~III.(1.1)]{hsiang},
so from the Serre spectral sequence 
of the bundle $\G/N_{\G}(T) \to B\G \to BN_{\G}(T)$
one sees $\H_\G \iso \H_{N_{\G}(T)}$.
From the principal bundle $N_{\G}(T)/T \to BT \to BN_{\G}(T)$,
the standard lemma on invariants 
under the action in cohomology induced by a covering 
action~\cite[Prop.~3G.1]{hatcher} 
shows $\HG \iso (\HT)^{N_\G(T)}$.

To apply this result in practice, 
we will need to understand the action of the normalizer on $\HT$. 
%It will help to see something more general.
First, we have an isomorphism $\Hom(T,S^1) \isoto H^1(T;\Z)$
given by pulling back the fundamental class in $H^1(S^1;\Z)$.
A homomorphism $\phi\: T \lt S^1$
induces a one-dimensional complex representation $V_\phi$,
and hence a line bundle $\defm{\xi_\phi}\: ET \ox_T V_\phi \lt BT$
with first Chern class $c_1(\xi_\phi) \in H^2(BT;\Z)$,
and $\phi \lmt c_1(\xi_\phi)$ is an isomorphism of abelian groups 
$\Hom(T,S^1) \isoto H^2(BT;\Z)$. 
The right covering action of an element of $N_\G(T)/T$ on $BT = E N_\G(T) / T$
is given by $eT \. nT = enT$,
and the group $\Hom(T,S^1)$ admits a right $N_\G(T)/T$-action
given by $\phi\. nT \: t \lmt \phi(n\-tn)$.
One easily checks  the pullback of 
$\xi_\phi$ along the right multiplication $eT \lmt enT$ on $BT$
is isomorphic to $\xi_{\phi\.nT}$,
so that the isomorphism 
$H^1(T;\Z) \isoto H^2(BT;\Z)$
is equivariant with respect to a left action of $N_\G(T)/T$.

Thus we can determine the action of $N_\G(T)/T$ on $\H_T$
by examining its conjugation action on $T$.
Moving back to our primary case of interest,
%
%
%Through this map, the normalizer $N_\G(T)$ of $T$
%in any supergroup of $T$ induces automorphisms of $\H(BT;\Z)$
%and we claim it is under this action that
%one can identify $\HK$ with an invariant subring of $\HT$.
%
%
%
%$BT \lt BK$
the map
$\HKz \longinc \HT$
is the tensor product of
the canonical injections $\H_{\SO(2k+\a)} \longmono \H_{T^k}$
and $\H_{\SO(2n+\b)} \longmono \H_{T^n}$.
The second factor can be viewed 
as the inclusion of invariants under the action of
either the Weyl group $\{\pm 1\}^{n} \semidirect S_n$ of $\SO(2n+1)$,
acting on $T^n$ by permutations and inversions of circle coordinates,
or the Weyl group of $\SO(2n)$, which is the subgroup 
	$\big\{(a_1,\ldots,a_n) \in \{\pm 1\}^n : 
	\prod a_j = 1\big\} \semidirect S_n$;
	and analogously for the first factor.
The induced action of the Weyl group on $\H(BT;\Z) \iso \Z[t_1,\ldots,t_n]$
is thus
by permutation and negation of the 
generators $\defm{t_j} = c_1(\xi_{\pr_j})$,
where $\smash{\pr_j\: T^n \lt S^1}$ is projection onto the $\smash{j\th}$ $\SO(2)$ factor. 
The invariants are generated by 
\eqn{
	p_\ell &= (-1)^\ell\s_\ell(t_1^2,\ldots,t_n^2),\\
	e  &= t_1 \cdots t_n,
}
where $\s_\ell$ is the $\ell\th$ elementary symmetric polynomial;
we can equivalently write $\smash{\sum p_\ell = \prod (1 - t\substack{2\\ \mn j\,})}$.
As suggested by the notation, these elements $p_\ell$
and $e$ can be viewed as the Pontrjagin and Euler classes.
The presence of the element $e$ 
in $\H B\SO(2n)$ but not in $\H B\SO(2n+1)$
is accounted for, from this point of view,
by the fact $t_1 \cdots t_n$ 
is invariant under the smaller Weyl group of $\SO(2n)$
but is sent to $-t_1 \cdots t_n$ by half the elements of 
the larger Weyl group of $\SO(2n+1)$.
In our special case, one does not need to localize all the way to $\Q$ coefficients to see $p$ and $e$ generate; in fact~\cite[Thm.~3.16]{VBKT},
\eqn{\label{eq:HBSO}
	\quotientmed
	{\H\gpp{B\SO(2n);\Z}\,}{\,2\mbox{-torsion}} &
	\,\iso\, 
	\Z[p_1,\ldots,p_{n-1},e],
	\phantom{{}_{n}}\quad \deg p_j = 4j, \ 
	\deg e = 2n,	\\
	\quotientmed
	{\H\gpp{B\SO(2n+1);\Z}\,}{\,2\mbox{-torsion}} &
	\,\iso\, 
	\Z[p_1,\ldots,p_{n-1},p_n],
	\quad \deg p_j = 4j.
}

The expressions for Pontrjagin and Euler classes in terms of maximal tori
makes it clear that if 
$\wt p$ denotes the total Pontrjagin class in $\H_{\SO(2k+2n+\a+\b)}$
and $\wt e$ the Euler class,
while $p,e$ and $p',e'$ continue to denote the corresponding
classes in $\H_{\SO(2k+\a)}$ and $\H_{\SO(2n+\b)}$
respectively,
then 
in the even case,
the inclusion $\SO(2k+\a) \x \SO(2n+\b) \longinc \SO(2k+2n+\a+\b)$ 
induces the restriction
\eqn{
	\H_{\SO(2k+2n+\a+\b)} &
	\longinc \H_{\SO(2k+\a) \x \SO(2n+\b)},\\
	\wt p = \prod_{1}^{n+k}(1-t_j^2) &\ \,=\ \,  \prod_{1}^{n}(1-t_j^2)\.\prod_{1}^{k}(1-t_{n+j}^2) 
	= pp',\\
	\wt e = \prod_1^{n+k} t_j &\ \,=\, \ \prod_1^{n} t_j \.\prod_1^{k} t_{n+j} = ee'.
}
of subrings of $\smash{\H_{T^k \x T^n}}$.
(The second restriction is just $0 = 0$ unless $\a = \b = 0$.)

In the odd case $\a = \b = 1$, the maximal torus of $\SO(2k+2n+2)$
is $\smash{\defm{\wt T} \iso T^k \x T^n \x S^1}$, 
where the last $S^1$ factor can be seen as the ``middle'' block-diagonal $\SO(2)$
that is not included in the block-diagonal $\SO(2k+1) \x \SO(1+2n)$.
The inclusion of $T \iso T^k \x T^n \x \{1\}$
then induces a map 
$\Z[t_1,\ldots,t_{k+n},t_{k+n+1}] = \H(B\wt T;\Z) \lt \H(BT;\Z)$
annihilating the last generator $t_{k+n+1}$
and preserving the other $t_j$.
It follows that in this case again $\wt p \lmt pp'$
under $\smash{\HG \lt \HKz}$
but now $\wt e = \smash{\prod_{j=1}^{n+k+1}} t_j \lmt \prod_{j=1}^{n+k} t_j \. 0 = 0$ despite the class $\wt e$ being itself nonzero.

Finally, to understand $\HK \lt \HKz$,
we need to examine the effect of conjugating $T$ 
by an element of $N_K(T)$ not lying in $N_{K_0}(T)$.
Such an element is a pair $(h,h') \in \mathrm S\big(\Orth(2k+\a) \x \Orth(2n+\b)\big) = K$ of square matrices such that $\det h = \det h' = -1$,
e.g., the pair with $h$ and $h'$ block-diagonal of the form
$\big[\begin{smallmatrix}0& 1\\ 1&0\end{smallmatrix}\big]
\+ [1] \+ \cdots \+ [1]$.
Conjugating $T^{k+n} = T^k \x T^n$ by this $(h,h')$
inverts one coordinate circle of $T^k$ and one of $T^n$,
so the induced action on $H_{T^{k+n}} = \Z[t_1,\ldots,t_{n+k}]$
sends the generator $t_1$ to $-t_1$ and 
the generator $t_{k+1}$ to $-t_{k+1}$,
fixing the other generators $t_j$.
This operation thus preserves $p$ and $p'$ 
but sends $e \mapsto -e$ and $e' \mapsto -e'$,
so $e$ and $e'$ 
are not invariant under the action of $N_K(T)$ 
but $ee'$ is (though it may also be zero).
Since $BK_0 \lt BK$ is a double covering, 
we have $\H\mn\big(BK;\Zt\big) \iso \H\mn\big(BK_0;\Zt\big){}^{\Z/2}$ 
%already  
%(i.e., inverting $2$ without passing to $\Q$ coefficients), 
so the image in $\H\mn\big(BT;\Zt\big)$ is 
\[
\case{
	\Zt\big[p,p',ee'\big] 
	& 
	\mbox{if }\a = \b = 0,
	\\
	\Zt[p,p'] 
	& 
	\mbox{otherwise,}
}
\]
where $e^2 = p_k$ and $(e')^2 = p'_n$
are absorbed notationally into the total classes $p$ and $p'$.
\end{discussion}

Understanding the map $\HG \lt \HKz$, 
we can prove He's theorem.

% (with rational coefficients) without explicit reference
%to the Kapovitch model we will introduce shortly.

%\bdefn
%One calls a continuous action of a topological group $\G$ on a space $X$
%\defd{equivariantly formal} if the map $\H_\G(X;\Q) \lt \H(X;\Q)$
%induced by the injection $X \lt E\G \ox_\G X\: x \lmt e_0 \ox x$
%for any fixed $e_0 \in E\G$ is surjective.
%\edefn

\begin{proof}[Proof of \Cref{thm:He}.]
	\hypertarget{pf1}{We} first recover the results with rational coefficients.
	
	For oriented Grassmannians in the even cases,
	we apply \Cref{thm:equalrank}
	to arrive at the general form 
	$\HKz \ox_\HG \HKz$.
	To see the presentation in terms of generators,
	note that
	$G = \SO(2k+2n+\a+\b)$ and $\G = H = K_0$ share a maximal torus,
	and recall from \Cref{rmk:restrictions} 
	that the Pontrjagin and Euler classes
	restrict as
	$\defm{\rho^*}\:\defm{\wt p} \lmt pp'$
	and 
	$\defm{\wt e} \lmt ee'$.
	Since $\wt p$ and $\wt e$ generate $\HG$,
	the tensor product can be computed explicitly
	by modding out of $\HKzt$ 
	the ideal generated by the homogeneous components of 
	$1 \ox \rho^*\wt p - \rho^*\wt p \ox 1$
	and %, in case $G = \SO(2k+2n)$, 
	the element $1 \ox \rho^* \wt e - \rho^* \wt e \ox 1$;
	which with the identifications 
	$p = 1 \ox p$ and $\pi = p \ox 1$, etc.,
	become $pp' -\pi\pi'$ and $ee'-\e\e'$.
%	but we again observe $\rho^*$ sends $\wt p \lmt pp'$
%	and %$\wt \pi \mapsto \pi\pi'$,\;
%	$\wt e \mapsto ee'$, %and $\wt\e \mapsto \e\e'$.
%	and so with the identifications $p = 1 \ox p$, $\pi = p \ox 1$, and so on,
%	we obtain the expressions for $\HKz(G/K_0)$ in the table.

	For oriented Grassmannians in the odd cases,
	we know the left action of $K_0$ on $G/K_0$ 
	is equivariantly formal~%
	\cite{goertsches2012isotropy,goertschesnoshari2016},
	so the result will follow from \Cref{thm:cohomisotropyformal} 
	once we see
	the image of $H^*(G/K_0) \lt H^* G$ is $\ext[\eta]$.
	By equivariant formality
	the images of $\H_{K_0}(G/K_0) \lt \H G$ and 
	$\H(G/K_0) \lt \H G$ agree,
	so it is enough to see the image of the latter is $\ext[\eta]$.%
	\footnote{\	Takeuchi already showed 
		this~\cite[p.~320]{takeuchi1962pontrjagin},
		but we want to recover his result as a corollary
		of \Cref{thm:He},
		so we derive it anew.
	}
	For this, we look at the Cartan model 
	$(\HKz \ox \H G, \bar d)$
	which is the case $U = 1 \x K_0$
	of the Kapovitch model.
	The differential $\bar d$ vanishes on $\HKz$
	and sends $1 \ox z \lmt 1 \ox \rho^* \tau z$
	for $z \in P \H G$,
%	The map 
%	first note the classifying map 
%	of the principal bundle 
%	$G \to G \ox_{K_0} EG \to BK_0 = K_0 \lq EG$
%	induces a bundle map from this bundle to the universal 
%	bundle $G \to EG \to BG$
%	and hence a map between the {\SSS}s
%	of these bundles.
	where the transgressions $\tau z$
	can be identified with the Pontrjagin and Euler classes of $\HG$.
	Per %and hence to nonzero elements of $\HKz$,
	\Cref{rmk:restrictions},
	the restrictions to $\HKz$ 
	of the Pontrjagin classes $\wt p_j$ are nonzero
	but the restriction of $\wt e$ is zero,
	meaning the element
	$\eta \in P \H G$ that transgresses to $\wt e$ 
	satisfies $\bar d (1 \ox \eta) = 0$,
	so 
	$\eta \in \H G$
	is the pullback of 
	$[1 \ox \eta] \in \H(\HKz \ox \H G,d) \iso \H(G/K_0)$
	along $G \lt G/K_0$
	and hence a generator of $\wP$.

	To obtain the expressions for the unoriented Grassmannians,
	note that $\KzGKz \lt \KzGK$ and
	$\KzGKz \lt \KGK$ are principal bundles with 
	respective fibers $1 \x \Z/2$ and $\Z/2 \x \Z/2$,
	so %\cite{eckmann1949coverings}
	\cite[Prop.~3G.1]{hatcher}
	we may identify $\HKz(G/K)$ and $\HK(G/K)$ 
	respectively with the $(1 \x \Z/2)$- and $(\Z/2 \x \Z/2)$-invariants
	of $\HKz (G/{K_0})$.
	But from \Cref{rmk:restrictions},
	the generator of the first $\Z/2$ factor 
	sends $\e \lmt -\e$ and $\e' \lmt -\e'$, 
	fixing $\pi,\pi',p,$ and $p'$,
	and likewise with the second $\Z/2$ factor and $e,e'$.

To refine the expressions to $\Zt$ coefficients,
%For the singular cohomology, recall that the Grassmannians can be realized
%as finite CW-subcomplexes of infinite Grassmannians
%$B\Orth(\ell)$ and $B\SO(\ell)$.
%The cohomology of the larger spaces contain 
%torsion only of order $2$ \cite[Thm.~3.2]{millerSO}%
%\cite[Thms.~A,B,12.1]{thomasBSO}\cite[Thms.~1.5,6]{brownBSO},
%so from the long exact sequence of a CW-pair,
%the same holds of the a finite Grassmannian $X$.
%Thus the coefficient homomorphism 
%$\H(X;\Z) / ($2$\mbox{-torsion}) \lt \H(X;\Q)$
%is injective,
%and as the Pontrjagin and Euler classes are integral,
%they generate the image.
recall that the torsion of the fiber and base of the 
bundle $G \to \KzGKz \to BK_0 \x BK_0$
is all of order two
\cite[Thm.~3.2]{millerSO}%
\cite[Thms.~A,~B,~12.1]{thomasBSO}\cite[Thms.~1.5, 1.6]{brownBSO}, 
so the \SSS tells us
the torsion of $\KzGKz$ must be of $2$-primary order 
as well~\cite[Lem.~1.9]{SSAT}.
As the covering maps 
$\KzGK \lt \KGK$ and
$\KzGKz \lt \KGK$ 
are respectively two- and four-sheeted,
they induce an injection 
in cohomology with $\Zt$ coefficients~\cite[Prop.~3G.1]{hatcher}.
Thus the maps $\H_\G\big(G/H;\Zt\big) \lt \H_\G(G/H;\Q)$
for 
$\G,H \in \{K,K_0\}$
are all injective,
and $\H_\G(G/H;\Q) \iso \H_\G\big(G/H;\Zt\big) \ox_{\,\Zt} \Q$
by the universal coefficient theorem.
As $\pi$, $\pi'$, $p$, $p'$, $\e$, $\e'$, $e$, $e'$
generate $\H_\G(G/H;\Q)$
and lie in the subring $\H_\G\big(G/H;\Zt\big)$,
it follows they generate the subring as well.
% $\H_\G\big(G/H;\Zt\big)$
%all lie in the domain and their images span the codomain.
%
%\tc{Run the Bockstein spectral sequence.}
\end{proof}

%\brmk
%The avoidance of the Kapovitch model in this proof is illusory,
%as the proof of \Cref{thm:cohomisotropyformal} uses this model.
%\ermk

\begin{history}\label{rmk:history}
	Leray already in 1949 had determined the cohomology 
	rings $\H\big(\wt G_\ell(\R^{\ell+m});\R\big)$ in the even cases%
	~\cite{lerayCR1949e}.
	Cartan found the Poincar\'{e} polynomials in both even and odd cases
	in the contribution to the Brussels 
	\emph{Colloque} on fiber spaces 
	where he introduces what was to become known 
	as the Cartan model of Borel equivariant cohomology~\cite[p.~70--1]{cartan1950transgression};
	it is clear from context that he also understood the 
	ring structures, but these were apparently less interesting 
	to him than the Betti numbers. 
	Borel's thesis, as already discussed, 
	identified a general technique
	for computation of the cohomology of total spaces of principal
	bundles and particularly of homogeneous spaces 
	not dependent on smoothness hypotheses or the Lie-theoretic
	considerations underlying the Cartan model,
	and as an example provided the cohomology rings 
	of the oriented real Grassmannians in the even case~\cite[p.~192]{borelthesis}.
	The first explicit statement of the ring structure
	in the odd case is due to Takeuchi~\cite{takeuchi1962pontrjagin},
	as an application of 
	\Cref{thm:formalCDGA}
	(connecting the ring structure of $\H(G/H)$
	with the regular sequence criterion;
	formality had not been defined yet);
%	describing
%	the real cohomology ring of a formal homogeneous space;
	this can be seen as an elaboration of the method from Borel's
	thesis, and is actually already to be found in Cartan~\cite[p.~70--1]{cartan1950transgression}
	as Takeuchi notes in a postscript 
	that he has belatedly been informed.
	
	The author is uncertain where the corresponding expressions 
	for the unoriented Grassmannians first appear,
	but they are easily calculated as invariant subrings.
	%of the oriented Grassmannians, as will be elaborated later,
	%and 
	This expression was conjectured in 2013 to hold by
	Casian and Kodama~\cite[Conjecture~6.1]{casiankodama2013}.
	Chen He found this conjecture and proved it, and the oriented counterpart,
	without recourse to the work discussed in the previous paragraph, 
	in 2016~\cite[5.2.1--2,\, 6.3.1]{he2016grassmannian}.
	Moreover, He derived these results as a consequence of the 
	rational Borel equivariant cohomology rings of $G_\ell(\R^{\ell+m})$
	and $\wt G_\ell(\R^{\ell+m})$
	with respect to the action of a maximal torus $T$
	in the stabilizer of a point under the canonical 
	action of $\SO(\ell+m)$,
	which had not appeared before.
	Independently, 
	Rustam Sadykov had found a very brief \emph{sui generis} geometric 
	proof of the rational non-equivariant cohomology of 
	$\wt G_\ell(\R^{\ell+m})$
	in 2015 \cite[Thm.~1]{sadykov2017} using the tautological 
	sphere bundles and the Gysin sequence, and in 2017
	expanded his paper to provide an elementary proof of He's result 
	upon the suggestion of a referee~\cite[Thm.~6]{sadykov2017}.
	
	He's preprint gives a ``Borel presentation'' for these rings 
	equivalent to that given in this note,
	which is in turn derived from a GKM-style presentation 
	relying on his own expansion of standard GKM-theory 
	to allow for not necessarily oriented and possibly odd-dimensional manifolds.
	This GKM-type presentation is rather involved and relies on
	a detailed understanding of the combinatorics of 
	the low-dimensional orbits of these actions. 
	The contribution of this section of the present note is simply to observe that the 
	``Borel presentation''
	is readily derivable from general principles
	(in terms of results
	%	is actually a ready consequence of already-available 
	%	structure theorems for 
	%	for the Borel equivariant cohomology
	%	$\HK(G/H;\Q)$ 
	%of a homogeneous space 
	which have been available 
	%for some time 
	in the even cases since 1993
	and
	in the odd case since 2004%
	~\cite{kumarvergne1993equivariant,kapovitch2002biquotients})
	and requires relatively little geometric information 
	about Grassmannians.
	Although the calculation of 
	$\H_{\mathrm{S}(\O(\ell) \x \O(m))} \big(G_{\ell}(\R^{\ell+m}) ;\Zt\big)$
	via the covering argument is trivial, it
	%		The calculation of for one of the 		
	%		equivariant cohomology of the oriented Grassmannians is, 
	%		although an easy application of the same results
	is original to this note as far as the author knows.	
	The main result of the note, \Cref{thm:cohomisotropyformal}, 
	was first proven in 2015
	in an unpublished post-defense revision of the author's dissertation,
	and has been mentioned without proof in one published article~\cite[Cor.~3.4]{carlsonfok2018}.
	A version of the present proof of He's \Cref{thm:He}
	has been on the arXiv since late 2016~\cite{carlson2016grassmannian}.
	
	The author would like to thank Paul Baum, 
	Chen He, 
	Vitali Kapovitch,
	and Rustam Sadykov
	for conversations and correspondence clarifying the chronology of their results.
\end{history}

\vspace{-1.5em}

% \cite{shult2010points}

%\printbibliography[heading=bibintoc]
{\footnotesize\bibliography{bibshort} }

\begin{thebibliography}{Onishchik}

\bibitem[Baum68]{baum1968homogeneous}
Paul~F. Baum.
\newblock {On the cohomology of homogeneous spaces}.
\newblock {\em Topology}, 7(1):15--38, 1968.
  \href {http://dx.doi.org/10.1016/0040-9383(86)90012-1}
 {\path{doi:10.1016/0040-9383(86)90012-1}}

\bibitem[Borel]{borelthesis}
Armand Borel.
\newblock {Sur la cohomologie des espaces fibr{\'e}s principaux et des espaces
  homog{\`e}nes de groupes de {L}ie compacts}.
\newblock {\em Ann. of Math. (2)}, 57(1):115--207, 1953.
\newblock  \url{http://web.math.rochester.edu/people/faculty/doug/otherpapers/Borel-Sur.pdf},
  \href {http://dx.doi.org/10.2307/1969728} {\path{doi:10.2307/1969728}}.

\bibitem[Bro82]{brownBSO}
Edgar~H. Brown, Jr.
\newblock The cohomology of {$BSO_n$} and {$BO_n$} with integer coefficients.
\newblock {\em Proc. Amer. Math. Soc.}, pages 283--288, 1982.
  \href {http://dx.doi.org/10.2307/2044298}
 {\path{doi:10.2307/2044298}}





\bibitem[Carl16]{carlson2016grassmannian}
Jeffrey~D. Carlson.
\newblock The {Borel} equivariant cohomology of real {Grassmannians}.
\newblock Nov 2016.
\newblock \href {http://arxiv.org/abs/1611.01175} {\path{arXiv:1611.01175}}.

\bibitem[Carl19]{carlson2018eqftorus}
Jeffrey~D. Carlson.
\newblock Equivariant formality of isotropic torus actions.
\newblock {\em J. Homotopy Relat. Struct.}, 14(1):199--234, 2019.
\newblock \href {http://arxiv.org/abs/1410.5740} {\path{arXiv:1410.5740}},
  \href {http://dx.doi.org/10.1007/s40062-018-0207-5}
  {\path{doi:10.1007/s40062-018-0207-5}}.

\bibitem[CarlF18]{carlsonfok2018}
Jeffrey~D. Carlson and {Chi}-{Kwong} Fok.
\newblock Equivariant formality of isotropy actions.
\newblock {\em J. Lond. Math. Soc.}, Mar 2018.
\newblock \href {http://arxiv.org/abs/1511.06228} {\path{arXiv:1511.06228}},
  \href {http://dx.doi.org/10.1112/jlms.12116} {\path{doi:10.1112/jlms.12116}}.

\bibitem[Cart51]{cartan1950transgression}
Henri Cartan.
\newblock {La transgression dans un groupe de {Lie} et dans un espace fibr{\'e}
  principal}.
\newblock In {\em {Colloque de topologie (espace fibr{\'e}s), {Bruxelles}
  1950}}, pages 57--71, Li{\`e}ge/Paris, 1951. Centre belge de recherches
  math{\'e}matiques, Georges Thone/Masson et companie.
\newblock \url{http://eudml.org/doc/112227}.


\bibitem[CK13]{casiankodama2013}
Luis Casian and Yuji Kodama.
\newblock On the cohomology of real {Grassmann} manifolds.
\newblock Sep 2013.
\newblock \href {http://arxiv.org/abs/1309.5520} {\path{arXiv:1309.5520}}.

\bibitem[Esch92]{eschenburg1992biquotient}
Jost-Hinrich Eschenburg.
\newblock {Cohomology of biquotients}.
\newblock {\em Manuscripta Math.}, 75(2):151--166, 1992.
\newblock \href {http://dx.doi.org/10.1007/BF02567078}
  {\path{doi:10.1007/BF02567078}}.

\bibitem[FHT]{FHT}
Yves F{\'e}lix, Steve Halperin, and Jean-Claude Thomas.
\newblock {\em {Rational homotopy theory}}, volume 205 of {\em {Grad. Texts in
  Math.}}
\newblock Springer, 2001.

\bibitem[FOT]{FOT}
Yves F{\'e}lix, John Oprea, and Daniel Tanr{\'e}.
\newblock {\em {Algebraic models in geometry}}, volume~17 of {\em {Oxford Grad.
  Texts Math.}}
\newblock Oxford Univ. Press, Oxford, 2008.
\newblock \url{www.maths.ed.ac.uk/~v1ranick/papers/tanre.pdf}.


\bibitem[Fra19]{franz2019homogeneous}
Matthias Franz.
\newblock The cohomology rings of homogeneous spaces.
\newblock Jul 2019.
\newblock \href {http://arxiv.org/abs/1907.04777} {\path{arXiv:1907.04777}}.



\bibitem[GKM98]{GKM1998}
Mark Goresky, Robert Kottwitz, and Robert MacPherson.
\newblock {Equivariant cohomology, {Koszul} duality, and the localization
  theorem}.
\newblock {\em Invent. Math.}, 131(1):25--83, 1998.
\newblock \url{http://math.ias.edu/~goresky/pdf/equivariant.jour.pdf},
  \href {http://dx.doi.org/10.1007/s002220050197}
  {\path{doi:10.1007/s002220050197}}.



\bibitem[G12]{goertsches2012isotropy}
Oliver Goertsches.
\newblock {The equivariant cohomology of isotropy actions on symmetric spaces}.
\newblock {\em Doc. Math.}, 17:79--94, 2012.
\newblock \url{emis.ams.org/journals/DMJDMV/vol-17/03.pdf}, \href
  {http://arxiv.org/abs/1009.4079} {\path{arXiv:1009.4079}}.


\bibitem[GN16]{goertschesnoshari2016}
Oliver Goertsches and Sam~Haghshenas Noshari.
\newblock Equivariant formality of isotropy actions on homogeneous spaces
  defined by Lie group automorphisms.
\newblock {\em J. Pure Appl. Algebra}, 220(5):2017--2028, 2016.
\newblock \href {http://arxiv.org/abs/1405.2655} {\path{arXiv:1405.2655}},
  \href {http://dx.doi.org/10.1016/j.jpaa.2015.10.013}
  {\path{doi:10.1016/j.jpaa.2015.10.013}}.

\bibitem[GuM]{gugenheimmay}
Victor~K.A.M. Gugenheim and J.~Peter May.
\newblock {\em On the theory and applications of differential torsion
  products}, volume 142 of {\em Mem. Amer. Math. Soc.}
\newblock Amer. Math. Soc., 1974.


\bibitem[Hatcher]{hatcher}
Allen Hatcher.
\newblock {\em {Algebraic topology}}.
\newblock Cambridge Univ. Press, 2002.
\newblock \url{http://math.cornell.edu/~hatcher/AT/ATpage.html}.


\bibitem[Hat04]{SSAT}
Allen Hatcher.
\newblock {\em {Spectral sequences in algebraic topology}}.
\newblock 2004
\newblock manuscript.
\newblock \url{http://math.cornell.edu/~hatcher/SSAT/SSATpage.html}.


\bibitem[Hat17]{VBKT}
Allen Hatcher.
\newblock {\em {Vector bundles and {K}-theory}}.
\newblock 2017
\newblock manuscript.
\newblock \url{http://math.cornell.edu/~hatcher/VBKT/VBpage.html}.

\bibitem[He16]{he2016grassmannian}
Chen He.
\newblock {{GKM} theory, characteristic classes and the equivariant cohomology
  ring of the real {Grassmannian}}.
\newblock Sep 2016.
\newblock \href {http://arxiv.org/abs/1609.06243} {\path{arXiv:1609.06243}}.

\bibitem[Hopf41]{hopf1941hopf}
Heinz Hopf.
\newblock {{\"U}ber eie {Topologie} der {Gruppen-Mannigfaltigkeiten} und ihre
  {Verallgemeinerungen}}.
\newblock {\em Ann. of Math.}, 42(1):22--52, 1941.
\newblock \url{http://jstor.org/stable/1968985}.

\bibitem[Hsiang]{hsiang}
Wu-Yi Hsiang.
\newblock {\em {Cohomology theory of topological transformation groups}}.
\newblock Springer, 1975.

\bibitem[Kap04]{kapovitch2002biquotients}
Vitali Kapovitch.
\newblock {A note on rational homotopy of biquotients}.
\newblock 2004.
\newblock \url{http://www.math.toronto.edu/vtk/biquotient.pdf}.

\bibitem[KV93]{kumarvergne1993equivariant}
Shrawan Kumar and Mich{\`e}le Vergne.
\newblock {Equivariant cohomology with generalized coefficients}.
\newblock {\em Ast{\'e}risque}, 215:109--204, 1993.
\newblock \url{http://numdam.org/item/AST_1993__215__109_0}

\bibitem[Ler49]{lerayCR1949e}
Jean Leray.
\newblock {D{\'e}termination, dans les cas non exceptionnels, de l'anneau de
  cohomologie de l'espace homog{\`e}ne quotient d'un groupe de {Lie} compact
  par un sous-groupe de m$\mathrm{\hat{e}}$me rang}.
\newblock {\em C. R. Acad. Sci. Paris}, 228(25):1902--1904, 1949.
\newblock \url{http://gallica.bnf.fr/ark:/12148/bpt6k31801/f1902}

\bibitem[Mill53]{millerSO}
Clair~E. Miller.
\newblock The topology of rotation groups.
\newblock {\em Ann. of Math.}, pages 90--114, 1953.
  \href {http://dx.doi.org/10.2307/1969727}
  {\path{doi:10.2307/1969727}}.



\bibitem[Munk74]{munkholm1974emss}
Hans~J. Munkholm.
\newblock {The {Eilenberg--Moore} spectral sequence and strongly homotopy
  multiplicative maps}.
\newblock {\em J. Pure Appl. Algebra}, 5(1):1--50, 1974.
  \href {http://dx.doi.org/10.1016/0022-4049(74)90002-4}
  {\path{doi:10.1016/0022-4049(74)90002-4}}.


\bibitem[Onischik]{onishchik}
Arkadi~L. Onishchik.
\newblock {\em {Topology of transitive transformation groups}}.
\newblock Johann Ambrosius Barth, 1994.

\bibitem[Sad17]{sadykov2017}
Rustam Sadykov.
\newblock Elementary calculation of the cohomology rings of real {Grassmann}
  manifolds.
\newblock {\em Pacific J. Math.}, 289(2):443--447, 2017.
\newblock \href {http://dx.doi.org/10.2140/pjm.2017.289.443}
  {\path{doi:10.2140/pjm.2017.289.443}}.




\bibitem[Shi96]{shiga1996equivariant}
Hiroo Shiga.
\newblock {Equivariant de {Rham} cohomology of homogeneous spaces}.
\newblock {\em J. Pure Appl. Algebra}, 106(2):173--183, 1996.
\newblock \href {http://dx.doi.org/10.1016/0022-4049(95)00018-6}
  {\path{doi:10.1016/0022-4049(95)00018-6}}.

\bibitem[Sing93]{singhof1993}
Wilhelm Singhof.
\newblock {On the topology of double coset manifolds}.
\newblock {\em Math. Ann.}, 297(1):133--146, 1993.
\newblock \href {http://dx.doi.org/10.1007/BF01459492}
  {\path{doi:10.1007/BF01459492}}

\bibitem[Tak62]{takeuchi1962pontrjagin}
Masaru Takeuchi.
\newblock {On {Pontrjagin} classes of compact symmetric spaces}.
\newblock {\em J. Fac. Sci. Univ. Tokyo Sect. I}, 9:313--328, 1962.
\newblock \url{http://dropbox.com/s/858u6qp89gmp88e/Pont_cpt_sym_sp%28Takeuchi%2C1962%29.pdf}

\bibitem[Tho60]{thomasBSO}
Emery Thomas.
\newblock On the cohomology of the real {Grassmann} complexes and the
  characteristic classes of $n$-plane bundles.
\newblock {\em Trans. Amer. Math. Soc.}, 96(1):67--89, 1960.
\newblock \href {http://dx.doi.org/10.2307/1993484 }
  {\path{doi:10.2307/1993484 }}



\bibitem[Tu10]{tu2010homogeneous}
Loring~W. Tu.
\newblock {Computing characteristic numbers using fixed points}.
\newblock In Peter~Robert Kotiuga, editor, {\em {A celebration of the
  mathematical legacy of {Raoul} {Bott}}}, volume~50 of {\em {C. R. M.
  Proceedings and Lecture Notes}}, pages 185--206. Centre de Recherches
  Math{\'e}matiques Montr{\'e}al, Amer. Math. Soc., 2010.
\newblock \href {http://arxiv.org/abs/math/0102013} {\path{arXiv:math/0102013}}.

\bibitem[Wolf77]{wolf1977homogeneous}
Joel~L. Wolf.
\newblock {The cohomology of homogeneous spaces}.
\newblock {\em Amer. J. Math.}, pages 312--340, 1977.

\end{thebibliography}

\nd\footnotesize{\textsc{
		%, 
		%Toronto, Ontario\ \ M5S 2E4, Canada
	}\\
	\url{jeffrey.carlson@tufts.edu}
	%\\\url{http://www.math.toronto.edu/jcarlson}
}

\end{document}